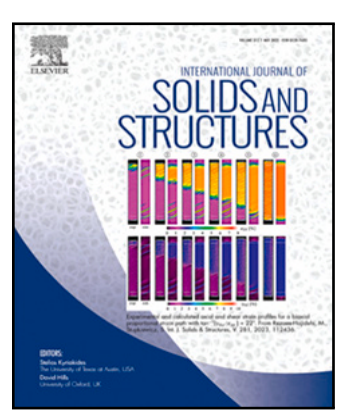

# Surface stability of a layered magnetoelastic half-space

D. Shahsavari [a], L. Dorfmann [b], P. Saxena [a],*

[a] Glasgow Computational Engineering Centre, James Watt School of Engineering, University of Glasgow, UK

[b] Department of Civil and Environmental Engineering, Tufts University, Medford, MA, USA

ABSTRACT

We evaluate the conditions for surface stability of a layered magnetoelastic half-space subjected to large deformations and a magnetic field. After reviewing the fundamental measures of deformation and summarizing the magnetostatic equations in Eulerian and Lagrangian forms, we derive the constitutive relations from a total energy function dependent on the deformation gradient and the Lagrangian magnetic induction. Energy principles yield the equilibrium equations, magnetic field equations, and boundary conditions. The second variation of the energy functional provides the incremental equations and conditions for stability analysis. Surface instability is studied by linearizing increments of deformation and magnetic induction about a finitely deformed state under a magnetic field normal to the surface. Four illustrative cases are considered: (i) a layered non-magnetizable half-space with varying stiffness contrast; (ii) the critical stretch of a magnetoelastic half-space as a function of magnetic induction; (iii) surface stability of a magneto-sensitive layer atop a non-magnetizable substrate; and (iv) bifurcation conditions of a two-layers magnetoelastic solid with different stiffness ratios. Bifurcation criteria in the form of coupled critical stretch and wavenumber is determined using a bespoke optimization protocol developed using an arc-length continuation method. Graphical results are provided throughout.

## 1. Introduction

The behavior of soft magnetizable materials has received increasing attention in recent years, as evidenced by their growing use in research and industrial applications. These materials, also known as magneto-active, magneto-sensitive or magneto-rheological elastomers (MRE), are utilized in soft robots, sensors, artificial muscles, and in biomedical devices (Chen et al., 2021; Zhao et al., 2021; Zhou et al., 2021; Silva et al., 2022; Bira et al., 2020; Kim and Zhao, 2022). The unique advantage is that the mechanical properties can be changed by an externally applied magnetic field.

The theory of magnetoelasticity, with constitutive equations specialized for the nonlinear response of magneto-sensitive materials, encompasses large deformations and magnetostatics. These two aspects have evolved independently, as their coupling has not been necessary.

The need to develop a large deformation theory and consistent stress–strain relations was primarily driven by the increased use of elastomeric materials in industrial applications. Significant contributions were made by Ronald Rivlin, who, in a series of articles, found new and elegant solutions to describe the behavior of elastomers and other soft solids (Rivlin, 1948a,b,c,d, 1949). In the resulting nonlinear theory, the field variables, equilibrium equations, and boundary conditions are defined in Eulerian and Lagrangian forms, with the corresponding vector transformations named after the Italian mathematician Gabrio Piola. A comprehensive summary of the fundamental concepts of rubber elasticity is given in Treloar (2005), a modern treatment of the mathematical theory and its application to the solution of boundary value problems, including bifurcation and stability in Ogden (1997).

The experimental works by Coulomb, Ampère and Faraday were used by Maxwell to formulate the theory on electricity and magnetism. Over a ten years period he presented his theoretical findings in a series of three papers (Maxwell, 1856, 1861, 1865). Maxwell's equations were derived in Cartesian coordinates and involve electric and magnetic fields and vector and scalar potentials. There are twenty of these equations in all, involving twenty variables. In 1884–1885 Oliver Heaviside rewrote Maxwell's equations using the operators grad, div and curl and, thereby, reduced the original set to the four we use today (Hampshire, 2018).

To accurately describe the intrinsic mechanical and magnetic coupling of magnetizable materials, the theory of magnetoelasticity must be formulated within the framework of nonlinear continuum mechanics. Significant contributions to the development of this theory are given in Brown (1966), Pao (1978), Maugin (1988) and Dorfmann and Ogden (2014). Recent reviews are given in Bira et al. (2020), Cantera

* Corresponding author.
*E-mail address:* Prashant.Saxena@glasgow.ac.uk (P. Saxena).

et al. (2017), Dorfmann and Ogden (2023) and Saber and Sedaghati (2023).

Constitutive equations that govern nonlinear magnetoelastic interactions are needed to illustrate the influence of the magnetic field on the mechanical response of magneto-sensitive elastomers. Notable contributions include those in Brigadnov and Dorfmann (2003), Steigmann (2004), Kankanala and Triantafyllidis (2004), Dorfmann and Ogden (2004b,a), Danas et al. (2012) and the implicit equations in Bustamante and Rajagopal (2014). A constitutive theory for transversely isotropic magnetoelastic materials is proposed in Bustamante (2010) and Saxena et al. (2014). Homogenization methods have been employed in Ponte Castañeda and Galipeau (2011), Chatzigeorgiou et al. (2014), Danas (2017) and Mukherjee et al. (2020). Numerical solutions of finite geometry boundary-value problems in nonlinear magnetoelasticity are given in Bustamante et al. (2007, 2011). The time dependent behavior of magneto-elastic material is considered in Saxena et al. (2013, 2014).

Stability analysis in magnetoelastic materials is typically reduced to studying linearized increments in the deformation and in the magnetic induction, superimposed on a finitely deformed configuration in the presence of a magnetic field. The formulation of Dorfmann and Ogden (2004a) is used in Otténio et al. (2008) as the starting point for the derivation of the incremental governing equations, the magnetoelastic moduli and the incremental boundary conditions. In Otténio et al. (2008) and Shahsavari and Saxena (2025) the general incremental equations are applied to the analysis of surface stability of a finitely deformed elastic half-space. In Saxena and Ogden (2011, 2012) the formulation is used to analyze surface and Love-type waves in a finitely deformed magnetoelastic half-space. The surface stability of a magneto-sensitive rectangular block and of a layer resting on a non-magnetizable substrate is considered, respectively, in Kankanala and Triantafyllidis (2008) and Danas and Triantafyllidis (2014). An experimental and numerical study of the stability and post-bifurcation of a magneto-sensitive film bonded to a passive elastomer substrate is provided in Psarra et al. (2019). A finite element analysis is used in Rambausek and Danas (2021) to investigate the loss of stability of single and multi-layered magneto-sensitive composites. Different topologies of magnetizable and elastomeric layers and different volume fractions of magnetic particles are considered.

A numerical analysis in Goshkoderia et al. (2020) shows that instabilities at different length-scales develop in a magneto-sensitive elastomer with periodically distributed magnetizable particles, which are not possible in the absence of a magnetic field. Instabilities of a composite consisting of a magneto-sensitive layer embedded into an elastomeric matrix in the presence of a magnetic field and subjected to compressive strains are investigated in Chen et al. (2023). The experiments and computations show that the buckling strain depends on the applied field.

We present a bifurcation analysis based on linearized incremental equations to detect the onset of instabilities. Creasing instability, a finite-amplitude localized event, requires a fully nonlinear analysis to track the finite-amplitude deformation path. We note that the variational framework and the null-space formulation provided constitute the foundation for post-bifurcation analysis, as indicated in Cao and Hutchinson (2012) for the purely mechanical half-space. The objective of this paper is to present a Lagrangian framework to analyze the stability of a magnetoelastic layered half-space. The incremental method developed in Otténio et al. (2008) and Shahsavari and Saxena (2025) is extended to analyze a layered half-space with magnetic load applied normal to the surface. A systematic numerical optimization approach is developed using an enriched arc-length method to simultaneously determine critical stretch and wave number as the bifurcation criteria.

Surface instabilities of purely elastic layered materials and of homogeneous magnetoelastic half-spaces have been extensively studied. In contrast, the case of layered magnetoelastic materials, where mechanical stiffness and magnetic coupling vary between layers, remains comparatively less explored. The present research fills this gap through the development of a unified Lagrangian bifurcation framework that supports coupled field equations and jump conditions at material interfaces, in addition to a continuation-based optimization approach that results in the critical stretch and wavenumber simultaneously. It allows for systematic research of the combined effects of stiffness ratio and magnetic characteristics on surface stability, including domains in which magnetic loading alone may result in surface patterning. The exact and asymptotic stability analyses available for coated or bilayer elastic half-spaces proved valuable in validating our numerical approach, see, for example, Alawiye et al. (2019), Cai and Fu (1999, 2000), Fu and Cai (2015), Fu and Ciarletta (2015) and Fu et al. (2021, 2024). Deriving asymptotic expressions for magnetoelastic problems provides significant mathematical challenges and is outside the scope of this work.

This paper is organized as follows. Section 2 summarizes the kinematic quantities and Section 3 the basic equations of magnetostatics. In Section 4, we introduce the constitutive equations for a nonlinear isotropic magnetoelastic material and define a modified strain energy function in terms of the Lagrangian magnetic induction and the deformation gradient. In Section 5 the first variation of an energy functional results in the equilibrium and magnetic field equations, along with the mechanical and magnetic boundary conditions. The second variation yields the incremental forms of the total nominal stress, the Lagrangian magnetic field, the governing equations, and boundary conditions. The boundary value problem is formulated in Section 6, including both the underlying deformation and the superimposed increments in deformation and magnetic induction. Section 7 outlines the solution procedure, and Section 8 presents the numerical results in graphical form. Finally, concluding remarks are provided in Section 9.

## 2. Kinematics

Consider a deformable magneto-sensitive body, which occupies the unstressed configuration $\mathcal{B}_{\mathrm{r}}$ with boundary $\partial\mathcal{B}_{\mathrm{r}}$, surrounded by free space $\mathcal{B}'_{\mathrm{r}}$ limited at infinity by the boundary $\partial\mathcal{B}^{\infty}$. The application of mechanical or magnetic forces induces a quasi-static deformation resulting in the current configuration $\mathcal{B}$, with boundary $\partial\mathcal{B}$ and exterior space $\mathcal{B}'$. We assume that the boundary at infinity is fixed. A material point in $\mathcal{B}_{\mathrm{r}}$ is identified by the position vectors $\mathbf{X}$, which in the current configuration occupies the position $\mathbf{x} = \chi(\mathbf{X})$, where the vector field $\chi$ describes the deformation. The corresponding gradient tensor $\mathbf{F}$ is calculated as

$$\mathbf{F} = \operatorname{Grad} \chi(\mathbf{X}), \tag{1}$$

where Grad is the gradient operator with respect to $\mathbf{X}$. We introduce the standard notation

$$J = \det \mathbf{F}, \tag{2}$$

which defines the ratio of an infinitesimal volume element in the current configuration to the corresponding volume in the reference configuration. Therefore, for an incompressible material $J = 1$. Associated with $\mathbf{F}$ we have the right and left Cauchy–Green deformation tensors, defined respectively by

$$\mathbf{C} = \mathbf{F}^{\mathrm{T}}\mathbf{F}, \quad \mathbf{b} = \mathbf{F}\mathbf{F}^{\mathrm{T}}. \tag{3}$$

For a summary of the main ingredients of the large deformation theory we refer to, for example, Dorfmann and Ogden (2014).

## 3. Basic equations of magnetostatics

We provide a brief overview of the equations of magnetostatics in Eulerian and Lagrangian forms and refer to Dorfmann and Ogden (2003, 2004a,b) and Sharma and Saxena (2020) for the derivation and to Dorfmann and Ogden (2023) for a review of the experimental and theoretical contributions. Solutions of some representative boundary value problems are given in Dorfmann and Ogden (2005).

In the current configuration $\mathcal{B}$ we denote by $\mathbf{B}$ and $\mathbf{H}$, respectively, the magnetic induction and the magnetic field. In the absence of distributed currents and time dependence these satisfy Maxwell's equations

$$\operatorname{div}\mathbf{B} = 0, \quad \operatorname{curl}\mathbf{H} = \mathbf{0}, \tag{4}$$

where curl and div are the curl and divergence operators with respect to $\mathbf{x}$. The fields $\mathbf{H}$ and $\mathbf{B}$ are connected by the magnetization $\mathbf{M}$ such that

$$\mathbf{B} = \mu_0\left(\mathbf{H} + \mathbf{M}\right), \tag{5}$$

where $\mu_0$ is the permeability of free space. In free space or in a non-magnetizable material Eq. (5) is replaced by

$$\mathbf{B}^{\star} = \mu_0\mathbf{H}^{\star}, \tag{6}$$

where the superscript $\star$ is used to identify the corresponding fields. The form of (4) is unchanged with $\mathbf{H}$ and $\mathbf{B}$ replaced by $\mathbf{H}^{\star}$ and $\mathbf{B}^{\star}$.

On the boundary $\partial\mathcal{B}$, when there is no free surface current, the fields $\mathbf{B}$ and $\mathbf{H}$ satisfy the jump conditions

$$\mathbf{n}\cdot[\![\mathbf{B}]\!] = 0, \quad \mathbf{n}\times[\![\mathbf{H}]\!] = \mathbf{0}, \tag{7}$$

where the open square brackets signify the difference in the enclosed quantity evaluated on the boundary $\partial\mathcal{B}$ and $\mathbf{n}$ the normal outward pointing unit vector.

To describe the nonlinear magneto-mechanical coupling in a large deformation setting we focus on the reference configuration $\mathcal{B}_{\mathrm{r}}$. The Lagrangian counterparts of the Eulerian fields $\mathbf{B}$ and $\mathbf{H}$ are denoted $\mathbf{B}_{\mathrm{L}}$ and $\mathbf{H}_{\mathrm{L}}$ and, following the derivations in Dorfmann and Ogden (2004b, 2014), are obtained as

$$\mathbf{B}_{\mathrm{L}} = J\mathbf{F}^{-1}\mathbf{B}, \quad \mathbf{H}_{\mathrm{L}} = \mathrm{F}^{\mathrm{T}}\mathbf{H}. \tag{8}$$

These satisfy the Lagrangian field equations

$$\operatorname{Div}\mathbf{B}_{\mathrm{L}} = 0, \quad \operatorname{Curl}\mathbf{H}_{\mathrm{L}} = \mathbf{0}, \tag{9}$$

where Div and Curl are the differential and curl operators with respect to $\mathbf{X}$. We also introduce the Lagrangian forms

$$\mathbf{B}^{\star}_{\mathrm{L}} = J\mathbf{F}^{-1}\mathbf{B}^{\star}, \quad \mathbf{H}^{\star}_{\mathrm{L}} = \mathbf{F}^{\mathrm{T}}\mathbf{H}^{\star}, \tag{10}$$

where, following Toupin (1956), we have introduced a fictional deformation field that is a smooth extension of the deformation from inside the body, i.e. $\mathbf{F}^{\star} = \mathbf{F}$. These satisfy the equations

$$\operatorname{Div}\mathbf{B}^{\star}_{\mathrm{L}} = 0, \quad \operatorname{Curl}\mathbf{H}^{\star}_{\mathrm{L}} = \mathbf{0} \tag{11}$$

and Eq. (6) becomes

$$\mathbf{B}^{\star}_{\mathrm{L}} = \mu_0 J\mathbf{C}^{-1}\mathbf{H}^{\star}_{\mathrm{L}}. \tag{12}$$

The Lagrangian form of the jump conditions (7) are

$$\mathbf{N}\cdot[\![\mathbf{B}_{\mathrm{L}}]\!] = 0, \quad \mathbf{N}\times[\![\mathbf{H}_{\mathrm{L}}]\!] = \mathbf{0}, \tag{13}$$

where $\mathbf{N}$ is the unit outward normal on $\partial\mathcal{B}_{\mathrm{r}}$.

## 4. Magnetoelastic constitutive equations

It is standard practice to define the properties of magneto-sensitive materials in terms of a total energy density function $\Omega$. It can be defined either in terms of the Lagrangian magnetic induction vector $\mathbf{B}_{\mathrm{L}}$ or, alternatively in terms of the Lagrangian magnetic field $\mathbf{H}_{\mathrm{L}}$ as the independent magnetic variable, coupled with the deformation gradient $\mathbf{F}$, see Dorfmann and Ogden (2004a,b, 2023) for details. For the present study, we focus on just one of these representations, namely $\Omega(\mathbf{F}, \mathbf{B}_{\mathrm{L}})$. The total nominal and the total Cauchy stresses for an unconstrained material are then calculated by

$$\mathbf{T} = \frac{\partial\Omega}{\partial\mathbf{F}}, \quad \boldsymbol{\tau} = J^{-1}\mathbf{F}\frac{\partial\Omega}{\partial\mathbf{F}}. \tag{14}$$

For an incompressible material Eq. (14) is augmented to give

$$\mathbf{T} = \frac{\partial\Omega}{\partial\mathbf{F}} - p\mathbf{F}^{-1}, \quad \boldsymbol{\tau} = \mathbf{F}\frac{\partial\Omega}{\partial\mathbf{F}} - p\mathbf{I}, \tag{15}$$

where $p$ is a Lagrange multiplier and $\mathbf{I}$ the identity tensor.

For unconstrained and incompressible materials, the Lagrangian and Eulerian forms of the magnetic field are equal and are given by

$$\mathbf{H}_{\mathrm{L}} = \frac{\partial\Omega}{\partial\mathbf{B}_{\mathrm{L}}}, \quad \mathbf{H} = \mathbf{F}^{-\mathrm{T}}\frac{\partial\Omega}{\partial\mathbf{B}_{\mathrm{L}}}. \tag{16}$$

### 4.1. Material symmetry

Suppose that mechanical forces deform the unloaded body from the reference configuration $\mathcal{B}_{\mathrm{r}}$ to the current configuration $\mathcal{B}$. Further, suppose that in the current configuration a magnetic induction $\mathbf{B}$ with reference counterpart $\mathbf{B}_{\mathrm{L}}$ is applied. This identifies a specific direction in $\mathcal{B}_{\mathrm{r}}$ and gives the material an anisotropic character similar to that of a transversely isotropic fiber-reinforced material. For discussions of material symmetry for elastic and magnetoelastic materials we refer to Spencer (1971) and Steigmann (2004), respectively.

If there is no other structure present, the material is said to be an isotropic magnetoelastic material and, by objectivity, $\Omega$ depends on $\mathbf{F}$ through $\mathbf{C}$ and on the structure tensor $\mathbf{B}_{\mathrm{L}}\otimes\mathbf{B}_{\mathrm{L}}$. We express this dependence in the form

$$\Omega = \Omega(\mathbf{C}, \mathbf{B}_{\mathrm{L}}\otimes\mathbf{B}_{\mathrm{L}}). \tag{17}$$

The energy function is then reduced to dependence on the principal invariants $I_1, I_2, I_3$ of $\mathbf{C}$ defined by

$$I_1 = \operatorname{tr}\mathbf{C}, \quad I_2 = \frac{1}{2}[I_1^2 - \operatorname{tr}(\mathbf{C}^2)], \quad I_3 = J^2, \tag{18}$$

where tr is the trace of a second-order tensor.

Following Spencer (1971) we note that the integrity basis for two symmetric second-order tensors in three dimensions includes the invariants of $\mathbf{C}$, together with invariants $I_4, I_5, I_6$, which depend on $\mathbf{B}_{\mathrm{L}}\otimes\mathbf{B}_{\mathrm{L}}$. These have the forms

$$I_4 = \mathbf{B}_{\mathrm{L}}\cdot\mathbf{B}_{\mathrm{L}}, \quad I_5 = \mathbf{B}_{\mathrm{L}}\cdot(\mathbf{C}\mathbf{B}_{\mathrm{L}}), \quad I_6 = \mathbf{B}_{\mathrm{L}}\cdot(\mathbf{C}^2\mathbf{B}_{\mathrm{L}}). \tag{19}$$

For an unconstrained material the explicit form of the total Cauchy stress $(14)_2$ is obtained as

$$\begin{aligned}\boldsymbol{\tau} = 2J^{-1}\big[&\Omega_1\mathbf{b} + \Omega_2(I_1\mathbf{b} - \mathbf{b}^2) + I_3\Omega_3\mathbf{I} + I_3\Omega_5\mathbf{B}\otimes\mathbf{B}\\ &+ I_3\Omega_6(\mathbf{b}\mathbf{B}\otimes\mathbf{B} + \mathbf{B}\otimes\mathbf{b}\mathbf{B})\big],\end{aligned} \tag{20}$$

where the notation $\Omega_i$, $i = 1, 2, \ldots, 6$ is used to signify differentiation with respect to $I_1, I_2, \ldots, I_6$ respectively.

Using $(16)_2$, a direct calculation leads to the Eulerian form of the magnetic field vector

$$\mathbf{H} = 2J\left(\Omega_4\mathbf{b}^{-1} + \Omega_5 + \Omega_6\mathbf{b}\right)\mathbf{B}. \tag{21}$$

For an incompressible material $J \equiv 1$ and (20) is replaced by

$$\boldsymbol{\tau} = 2\Omega_1\mathbf{b} + 2\Omega_2(I_1\mathbf{b} - \mathbf{b}^2) - p\mathbf{I} + 2\Omega_5\mathbf{B}\otimes\mathbf{B} + 2\Omega_6(\mathbf{b}\mathbf{B}\otimes\mathbf{B} + \mathbf{B}\otimes\mathbf{b}\mathbf{B}), \tag{22}$$

which can also be expressed as the nominal stress

$$\begin{aligned}\mathbf{T} = {}& 2\Omega_1\mathbf{F}^{\mathrm{T}} + 2\Omega_2(I_1\mathbf{F}^{\mathrm{T}} - \mathbf{F}^{\mathrm{T}}\mathbf{F}\mathbf{F}^{\mathrm{T}}) - p\mathbf{F}^{-1}\\ &+ 2\Omega_5\mathbf{B}_{\mathrm{L}}\otimes\mathbf{F}\mathbf{B}_{\mathrm{L}} + 2\Omega_6(\mathbf{C}\mathbf{B}_{\mathrm{L}}\otimes\mathbf{F}\mathbf{B}_{\mathrm{L}} + \mathbf{B}_{\mathrm{L}}\otimes\mathbf{F}\mathbf{C}\mathbf{B}_{\mathrm{L}}).\end{aligned} \tag{23}$$

Similarly, from (21) we find that

$$\mathbf{H}_{\mathrm{L}} = 2\left(\Omega_4 + \Omega_5\mathbf{C} + \Omega_6\mathbf{C}^2\right)\mathbf{B}_{\mathrm{L}}. \tag{24}$$

## 5. Energy consideration

The objective is to construct a variational principle that produces the mechanical equations of equilibrium and boundary conditions together with the appropriate magnetic field and continuity equations.

We consider the energy functional in Lagrangian form and, following the development in Bustamante et al. (2008), Bustamante and Ogden (2012) and Sharma and Saxena (2020), we write

$$\begin{aligned} \Pi\left(\mathbf{x},\mathbf{A}_{\mathrm{L}}\right) = & \int_{\mathcal{B}_{\mathrm{r}}}\left[\Omega\left(\mathbf{F},\mathbf{B}_{\mathrm{L}}\right)-p(J-1)\right]\mathrm{d}V+\int_{\mathcal{B}_{\mathrm{r}}'}\Omega_{\mathrm{e}}\left(\mathbf{F},\mathbf{B}_{\mathrm{L}}^{\star}\right)\mathrm{d}V \\ & -\int_{\partial\mathcal{B}_{\mathrm{r}}}\mathbf{t}_{\mathrm{A}}\cdot\mathbf{x}\,\mathrm{d}A-\int_{\mathcal{B}_{\mathrm{r}}}\rho_{\mathrm{r}}\mathbf{f}\cdot\mathbf{x}\,\mathrm{d}V+\int_{\partial\mathcal{B}^{\infty}}\left(\mathbf{H}_{\mathrm{a}}\times\mathbf{A}\right)\cdot\mathbf{n}\,\mathrm{d}a, \end{aligned} \tag{25}$$

where the deformation $\mathbf{x}$ and the vector potential $\mathbf{A}_{\mathrm{L}}$ are the independent variables, $\mathbf{t}_{\mathrm{A}}$ is the applied mechanical traction per unit reference area and $\mathbf{f}$ the body force per unit mass. The magnetic induction is calculated as $\mathbf{B}_{\mathrm{L}} = \mathrm{Curl}\,\mathbf{A}_{\mathrm{L}}$, which automatically satisfies (9)$_1$. As in Bustamante et al. (2008), we assume that the vector potential is continuous across the boundary $\partial\mathcal{B}_{\mathrm{r}}$ resulting in $\mathbf{B}_{\mathrm{L}}^{\star} = \mathrm{Curl}\,\mathbf{A}_{\mathrm{L}}^{\star}$ and in the jump condition (13)$_1$ being satisfied. The Eulerian and Lagrangian forms of the vector potential are connected by $\mathbf{A} = \mathbf{F}^{-\mathrm{T}}\mathbf{A}_{\mathrm{L}}$, leading to $\mathbf{B} = \mathrm{curl}\,\mathbf{A}$ and (4)$_1$ is satisfied.

The energy $\Omega_{\mathrm{e}}\left(\mathbf{F},\mathbf{B}_{\mathrm{L}}^{\star}\right)$ in the exterior space $\mathcal{B}_{\mathrm{r}}'$ has the form

$$\Omega_{\mathrm{e}}\left(\mathbf{F},\mathbf{B}_{\mathrm{L}}^{\star}\right)=\frac{1}{2}\mu_0^{-1}J^{-1}\left(\mathbf{F}\mathbf{B}_{\mathrm{L}}^{\star}\right)\cdot\left(\mathbf{F}\mathbf{B}_{\mathrm{L}}^{\star}\right), \tag{26}$$

where $\mathbf{F} = \mathbf{F}^{\star}$ is the fictional deformation field introduced in Toupin (1956).

The prescribed field $\mathbf{H}_{\mathrm{a}} = \mu_0^{-1}\mathbf{B}_{\mathrm{a}}$ on the boundary $\partial\mathcal{B}^{\infty}$ is left in Eulerian form since the boundary at infinity is fixed.

### 5.1. First variation

We find that the functional $\Pi$ is stationary with respect to the independent variations $\mathbf{x}$ and $\mathbf{A}_{\mathrm{L}}$ if and only if

$$\mathrm{Div}\,\mathbf{T}+\rho_{\mathrm{r}}\mathbf{f}=\mathbf{0},\quad \mathrm{Curl}\,\mathbf{H}_{\mathrm{L}}=\mathbf{0},\quad J=1\quad\text{in}\quad\mathcal{B}_{\mathrm{r}}, \tag{27}$$

where definition (1) and the Lagrangian forms of the constitutive Eqs. (15)$_1$ and (16)$_1$ have been used. To find the variation in $\mathcal{B}_{\mathrm{r}}'$ it is useful to obtain the expression for

$$\begin{aligned} \frac{\partial\Omega_{\mathrm{e}}}{\partial\mathbf{F}} &= \mu_0^{-1}J^{-1}\mathbf{B}_{\mathrm{L}}^{\star}\otimes\mathbf{F}\mathbf{B}_{\mathrm{L}}^{\star}-\frac{1}{2}\mu_0^{-1}J^{-1}\left(\mathbf{F}\mathbf{B}_{\mathrm{L}}^{\star}\right)\cdot\left(\mathbf{F}\mathbf{B}_{\mathrm{L}}^{\star}\right)\mathbf{F}^{-1}, \\ &= \mu_0^{-1}J\mathbf{F}^{-1}\left[\mathbf{B}^{\star}\otimes\mathbf{B}^{\star}-\frac{1}{2}\left(\mathbf{B}^{\star}\cdot\mathbf{B}^{\star}\right)\mathbf{I}\right]=J\mathbf{F}^{-1}\boldsymbol{\tau}^{\star}, \end{aligned} \tag{28}$$

where $\boldsymbol{\tau}^{\star}$ is the Maxwell stress outside the material. It is convenient to introduce the Lagrangian form, denoted $\mathbf{T}^{\star}$, by

$$\mathbf{T}^{\star}=J\mathbf{F}^{-1}\boldsymbol{\tau}^{\star}, \tag{29}$$

which is similar to the transformation (14), which is valid for the purely mechanical case. In addition, we also calculate

$$\frac{\partial\Omega_{\mathrm{e}}}{\partial\mathbf{B}_{\mathrm{L}}^{\star}}=\mu_0^{-1}J^{-1}\mathbf{F}^{\mathrm{T}}\mathbf{F}\mathbf{B}_{\mathrm{L}}^{\star}=\mathbf{H}_{\mathrm{L}}^{\star}, \tag{30}$$

which coincides with (12).

In the exterior space we find that

$$\mathrm{Div}\,\mathbf{T}^{\star}=\mathbf{0},\quad \mathrm{Curl}\,\mathbf{H}_{\mathrm{L}}^{\star}=\mathbf{0}\quad\text{in}\quad\mathcal{B}_{\mathrm{r}}'. \tag{31}$$

We also obtain the mechanical and magnetic boundary conditions

$$\left(\mathbf{T}^{\star}-\mathbf{T}^{\mathrm{T}}\right)\mathbf{N}=\mathbf{t}_{\mathrm{A}}\quad\text{on}\quad\partial\mathcal{B}_{\mathrm{r}}, \tag{32}$$

$$\left(\mathbf{H}_{\mathrm{L}}^{\star}-\mathbf{H}_{\mathrm{L}}\right)\times\mathbf{N}=\mathbf{0}\quad\text{on}\quad\partial\mathcal{B}_{\mathrm{r}}, \tag{33}$$

and on the boundary $\partial\mathcal{B}^{\infty}$ we have

$$\left(\mathbf{H}_{\mathrm{a}}-\mathbf{H}^{\star}\right)\times\mathbf{n}=\mathbf{0}\quad\text{on}\quad\partial\mathcal{B}^{\infty}. \tag{34}$$

The boundary conditions (7)$_1$ and (13)$_1$ are satisfied automatically.

### 5.2. Second variation

We use the first variation results and begin with the second variation in $\mathbf{x}$. This gives

$$\int_{\mathcal{B}_{\mathrm{r}}}\mathrm{Div}\,\dot{\mathbf{T}}\,\mathrm{d}V-\int_{\partial\mathcal{B}_{\mathrm{r}}}\dot{\mathbf{T}}^{\mathrm{T}}\mathbf{N}\,\mathrm{d}A=\mathbf{0}, \tag{35}$$

where the increment of the nominal stress $\dot{\mathbf{T}}$ has the explicit form

$$\dot{\mathbf{T}}=\boldsymbol{\mathcal{A}}\dot{\mathbf{F}}+\boldsymbol{\Gamma}\dot{\mathbf{B}}_{\mathrm{L}}-\dot{p}\mathbf{F}^{-1}+p\mathbf{F}^{-1}\dot{\mathbf{F}}\mathbf{F}^{-1}, \tag{36}$$

the increment of the deformation gradient is denoted $\dot{\mathbf{F}}$ and the increment of the magnetic induction

$$\dot{\mathbf{B}}_{\mathrm{L}}=\mathrm{Curl}\,\dot{\mathbf{A}}_{\mathrm{L}}, \tag{37}$$

which leads to

$$\mathrm{Div}\,\dot{\mathbf{B}}_{\mathrm{L}}=0\quad\text{in}\quad\mathcal{B}_{\mathrm{r}}. \tag{38}$$

For reference we give $\dot{\mathbf{T}}$ in component form

$$\dot{T}_{\alpha i}=\mathcal{A}_{\alpha i\beta j}\dot{F}_{j\beta}+\Gamma_{\alpha i\beta}\dot{B}_{\mathrm{L}\beta}-\dot{p}F_{\alpha i}^{-1}+pF_{\alpha j}^{-1}\dot{F}_{j\beta}F_{\beta i}^{-1}, \tag{39}$$

where

$$\mathcal{A}_{\alpha i\beta j}=\frac{\partial^2\Omega}{\partial F_{i\alpha}\partial F_{j\beta}},\qquad \Gamma_{\alpha i\beta}=\frac{\partial^2\Omega}{\partial F_{i\alpha}\partial B_{\mathrm{L}\beta}},\qquad \dot{B}_{\mathrm{L}\beta}=\varepsilon_{\beta\gamma\alpha}\frac{\partial\dot{A}_{\alpha}}{\partial X_{\gamma}}. \tag{40}$$

To satisfy (35) both terms must vanish, for the first we have

$$\mathrm{Div}\,\dot{\mathbf{T}}=\mathbf{0}\quad\text{in}\quad\mathcal{B}_{\mathrm{r}}. \tag{41}$$

From the incompressibility constraint $J = 1$ follows that

$$\dot{J}=\mathbf{F}^{-1}\dot{\mathbf{F}}=0. \tag{42}$$

To find the incremental forms of the magnetic field and of the boundary condition we use the transformation

$$\int_{\mathcal{B}_{\mathrm{r}}}\dot{\mathbf{H}}_{\mathrm{L}}\cdot\dot{\mathbf{B}}_{\mathrm{L}}\,\mathrm{d}V=\int_{\mathcal{B}_{\mathrm{r}}}\mathrm{Curl}\,\dot{\mathbf{H}}_{\mathrm{L}}\cdot\dot{\mathbf{A}}_{\mathrm{L}}+\int_{\partial\mathcal{B}_{\mathrm{r}}}\left(\dot{\mathbf{H}}_{\mathrm{L}}\times\mathbf{N}\right)\cdot\dot{\mathbf{A}}_{\mathrm{L}}\,\mathrm{d}A, \tag{43}$$

where

$$\dot{\mathbf{H}}_{\mathrm{L}}=\boldsymbol{\Gamma}^{\mathrm{T}}\dot{\mathbf{F}}+\boldsymbol{\mathcal{K}}\dot{\mathbf{B}}_{\mathrm{L}}. \tag{44}$$

For reference, we write (44) in component form

$$\dot{H}_{\mathrm{L}\beta}=\Gamma_{\alpha i\beta}\dot{F}_{i\alpha}+\mathcal{K}_{\beta\alpha}\dot{B}_{\mathrm{L}\alpha}, \tag{45}$$

where the second order tensor $\boldsymbol{\mathcal{K}}$ is calculated as

$$\mathcal{K}_{\beta\alpha}=\frac{\partial^2\Omega}{\partial B_{\mathrm{L}\beta}\partial B_{\mathrm{L}\alpha}}. \tag{46}$$

We require that both terms on the right hand side of (43) vanish, from the first follows that

$$\mathrm{Curl}\,\dot{\mathbf{H}}_{\mathrm{L}}=\mathbf{0}\quad\text{in}\quad\mathcal{B}_{\mathrm{r}}. \tag{47}$$

In the surrounding space following an identical process we have

$$\mathrm{Div}\,\dot{\mathbf{T}}^{\star}=\mathbf{0}\quad\text{in}\quad\mathcal{B}_{\mathrm{r}}', \tag{48}$$

where the Lagrangian form of the Maxwell stress increment is obtained as

$$\dot{\mathbf{T}}^{\star}=\mathbf{F}^{-1}\dot{\boldsymbol{\tau}}^{\star}-\mathbf{F}^{-1}\dot{\mathbf{F}}^{\star}\mathbf{F}^{-1}\boldsymbol{\tau}^{\star}, \tag{49}$$

combined with

$$\dot{\boldsymbol{\tau}}^{\star}=\mu_0^{-1}\left[\dot{\mathbf{B}}^{\star}\otimes\mathbf{B}^{\star}+\mathbf{B}^{\star}\otimes\dot{\mathbf{B}}^{\star}-\left(\mathbf{B}^{\star}\cdot\dot{\mathbf{B}}^{\star}\right)\mathbf{I}\right]. \tag{50}$$

In addition, using (30) we obtain

$$\dot{\mathbf{H}}_{\mathrm{L}}^{\star}=\mu_0^{-1}\left(\dot{\mathbf{C}}\mathbf{B}_{\mathrm{L}}^{\star}+\mathbf{C}\dot{\mathbf{B}}_{\mathrm{L}}^{\star}\right), \tag{51}$$

where

$$\dot{\mathbf{B}}_{\mathrm{L}}^{\star}=\mathbf{F}^{-1}\dot{\mathbf{B}}^{\star}-\mathbf{F}^{-1}\dot{\mathbf{F}}^{\star}\mathbf{F}^{-1}\mathbf{B}^{\star} \tag{52}$$

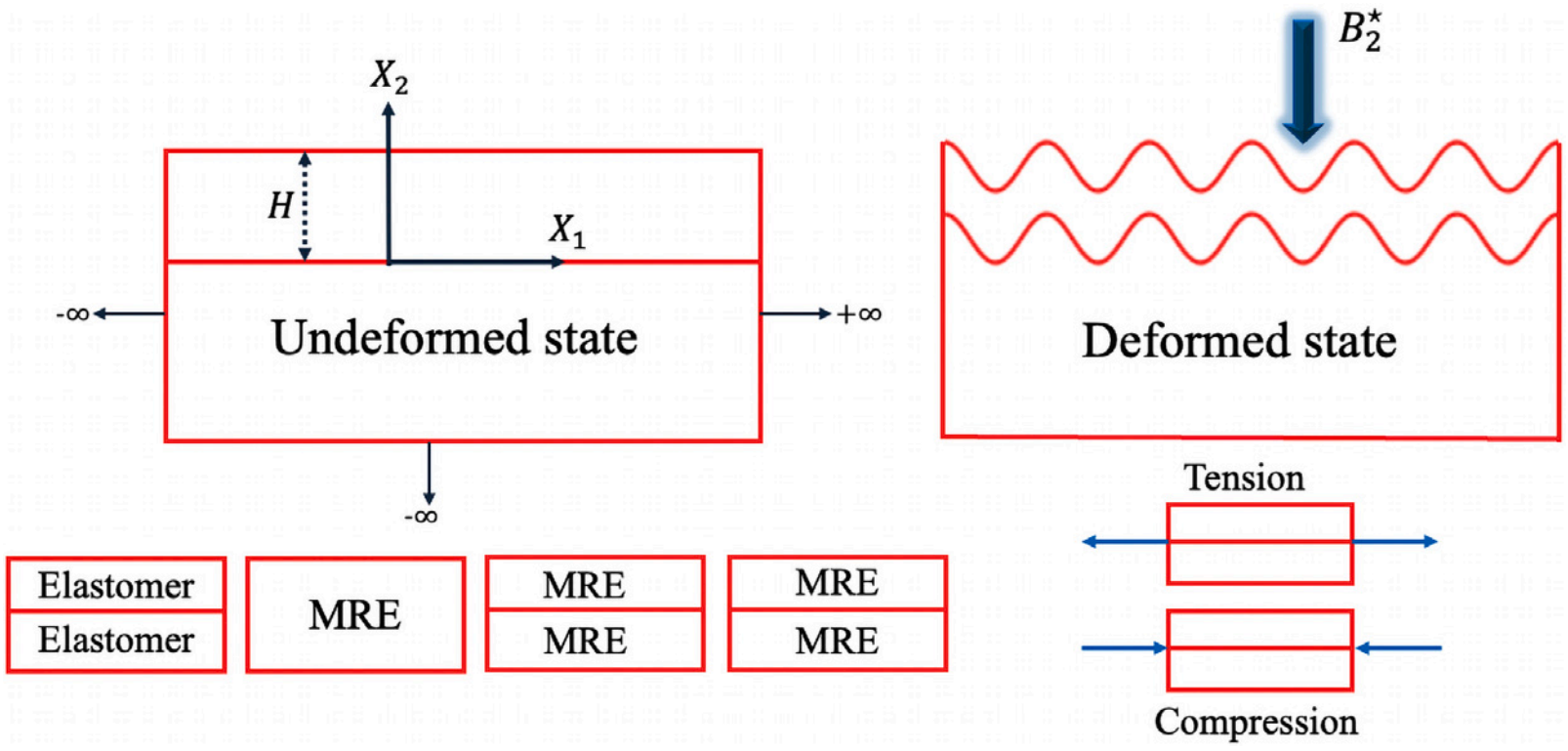

**Fig. 1.** A two-layer half-space is subjected to in-plane compression or tension, with a uniform magnetic induction normal to its boundary surface prior to the application of superimposed increments. First, we examine a two-layer non-magnetizable elastomeric solid. Next, we analyze a magnetoelastic half-space. Third, we investigate a non-magnetizable elastic substrate combined with a magnetizable upper layer. Finally, we study a layered magnetoelastic assembly.

and

$$\dot{\mathbf{H}}^\star_{\mathrm{L}} = \mathbf{F}^{\mathrm{T}}\dot{\mathbf{H}}^\star + \dot{\mathbf{F}}^{\star\mathrm{T}}\mathbf{H}^\star. \tag{53}$$

These Lagrangian fields satisfy the field equations

$$\mathrm{Div}\,\dot{\mathbf{B}}^\star_{\mathrm{L}} = 0, \quad \mathrm{Curl}\,\dot{\mathbf{H}}^\star_{\mathrm{L}} = \mathbf{0}. \tag{54}$$

We assume that the increment of the deformation gradient in free space differs from the increment in the material, hence the notation $\dot{\mathbf{F}}^\star$ is used to differentiate from $\dot{\mathbf{F}}$. It follows that $\dot{\mathbf{C}} = \dot{\mathbf{F}}^{\star\mathrm{T}}\mathbf{F} + \mathbf{F}^{\mathrm{T}}\dot{\mathbf{F}}^\star$.

The vanishing of the second term in (35) combined with the equivalent expression for the surrounding space gives the incremental traction boundary condition

$$\left(\dot{\boldsymbol{\tau}}^\star\mathbf{F}^{-\mathrm{T}} + \boldsymbol{\tau}^\star\mathbf{F}^{-\mathrm{T}}\dot{\mathbf{F}}^{\star\mathrm{T}}\mathbf{F}^{-\mathrm{T}} - \dot{\mathbf{T}}^{\mathrm{T}}\right)\mathbf{N} = \dot{\mathbf{t}}_{\mathrm{A}}. \tag{55}$$

Using (43) combined with (53) results in

$$\left(\mathbf{F}^{\mathrm{T}}\dot{\mathbf{H}}^\star + \dot{\mathbf{F}}^{\star\mathrm{T}}\mathbf{H}^\star - \dot{\mathbf{H}}_{\mathrm{L}}\right) \times \mathbf{N} = \mathbf{0}. \tag{56}$$

The continuity of the increments of the vector potential gives

$$\left(\mathbf{F}^{-1}\dot{\mathbf{B}}^\star - \mathbf{F}^{-1}\dot{\mathbf{F}}^\star\mathbf{F}^{-1}\mathbf{B}^\star - \dot{\mathbf{B}}_{\mathrm{L}}\right) \cdot \mathbf{N} = 0. \tag{57}$$

## 6. Problem description

In this section we focus on applying the foregoing theory to the pure homogeneous deformation of a half-space subjected to a magnetic field oriented perpendicular to the boundary surface. To investigate surface stability, we use the equations governing the linearized increments in deformation and in magnetic induction, along with the associated mechanical and magnetic boundary conditions (Fig. 1). The surface stability of a magnetizable half-space is considered in Otténio et al. (2008) and Shahsavari and Saxena (2025) and that of a magnetizable upper layer combined with a non-magnetic substrate is presented in Danas and Triantafyllidis (2014).

We use the rectangular Cartesian coordinates $X_1, X_2, X_3$ to define the reference configuration $\mathcal{B}_{\mathrm{r}}$ by

$$-\infty < X_1 < \infty, \quad -\infty < X_2 \leq H, \quad -\infty < X_3 < \infty, \tag{58}$$

where $H$ denotes the undeformed thickness of the upper layer occupying the region $0 \leq X_2 \leq H$. We take the material to be incompressible and subject to plane strain in the $(X_1, X_2)$ plane.

Consider the underlying homogeneous deformation defined by

$$x_1 = \lambda_1 X_1, \quad x_2 = \lambda_2 X_2, \quad x_3 = X_3. \tag{59}$$

where the incompressibility condition $\lambda_1\lambda_2 = 1$ suggests to introduce the notation $\lambda_1 = \lambda, \lambda_2 = \lambda^{-1}$. The deformation gradient $\mathbf{F}$ is then given by the matrix of Cartesian components $\mathsf{F}$ as

$$\mathsf{F} = \begin{bmatrix} \lambda & 0 & 0 \\ 0 & \lambda^{-1} & 0 \\ 0 & 0 & 1 \end{bmatrix} \tag{60}$$

and, with the use of (18), we obtain the invariants

$$I_1 = I_2 = \lambda^2 + \lambda^{-2} + 1, \quad I_3 = 1. \tag{61}$$

We start by applying a magnetic induction vector $\mathbf{B}_{\mathrm{L}}$ in the $X_2$ direction with component $B_{\mathrm{L}2}$. From Eq. (9) and from the boundary condition $(13)_1$ we find that $B_{\mathrm{L}2}$ is constant and equal in the substrate, in the upper layer and in the surrounding space

$$B_{\mathrm{L}2\mathrm{s}} = B_{\mathrm{L}2\mathrm{u}} = B^\star_{\mathrm{L}2}, \tag{62}$$

where the subscripts s and u indicate the substrate and the upper layer, respectively. This results in the invariants (19) having the forms

$$I_4 = B^2_{\mathrm{L}2}, \quad I_5 = \lambda^{-2}I_4, \quad I_6 = \lambda^{-4}I_4. \tag{63}$$

Using $B^\star_2 = \lambda^{-1}B^\star_{\mathrm{L}2}$ in (28) gives the Maxwell stress components exterior to the material. These are

$$\tau^\star_{22} = \frac{1}{2}\mu_0^{-1}\lambda^{-2}I_4, \quad \tau^\star_{11} = \tau^\star_{33} = -\frac{1}{2}\mu_0^{-1}\lambda^{-2}I_4. \tag{64}$$

The corresponding Lagrangian forms are obtained using (29) and result in

$$T^\star_{11} = -\frac{1}{2}\mu_0^{-1}\lambda^{-3}I_4, \quad T^\star_{22} = \frac{1}{2}\mu_0^{-1}\lambda^{-1}I_4, \quad T^\star_{33} = -\frac{1}{2}\mu_0^{-1}\lambda^{-2}I_4. \tag{65}$$

We now specialize the stress field to the plane strain condition using (23), (60) and (62). The resulting nonzero components of the nominal stress are

$$\begin{aligned} T_{11} &= 2\Omega_1\lambda + 2\Omega_2(\lambda + \lambda^{-1}) - p\lambda^{-1}, \\ T_{22} &= 2\Omega_1\lambda^{-1} + 2\Omega_2(\lambda + \lambda^{-1}) - p\lambda + 2\Omega_5\lambda^{-1}I_4 + 4\Omega_6\lambda^{-3}I_4, \\ T_{33} &= 2\Omega_1 + 2\Omega_2(\lambda^2 + \lambda^{-2}) - p. \end{aligned} \tag{66}$$

The traction continuity condition requires that the component $T_{22}$ evaluate in the two substrates is continuous across the material interface. In addition, with no applied mechanical traction, the boundary condition (32) requires that that component $T_{22}$ in the material is equal to the Maxwell component $(65)_2$. Once the energy function is specified, these two conditions will be used to evaluate the Lagrangian multipliers in the layered material.

Using (24) we obtain the nonzero component of the magnetic field as

$$H_{\mathrm{L}2} = 2\left(\Omega_4 + \lambda^{-2}\Omega_5 + \lambda^{-4}\Omega_6\right)B_{\mathrm{L}2}. \tag{67}$$

To analyze surface stability we define a displacement increment $\dot{\mathbf{u}}$ in the $(X_1, X_2)$ plane such that the incremental deformation gradient becomes

$$\dot{\mathbf{F}} = \operatorname{Grad}\dot{\mathbf{u}} = \begin{bmatrix} \dot{u}_{1,1} & \dot{u}_{1,2} & 0 \\ \dot{u}_{2,1} & \dot{u}_{2,2} & 0 \\ 0 & 0 & 0 \end{bmatrix}, \tag{68}$$

where $\dot{u}_{i,\alpha} = \partial \dot{u}_i / \partial X_\alpha$. It follows that the incompressibility condition (42) has the explicit form

$$\dot{u}_{1,1} + \lambda^2 \dot{u}_{2,2} = 0. \tag{69}$$

In addition to the incremental deformation, we consider an increment in the magnetic induction $\dot{\mathbf{B}}_{\mathrm{L}}$ with components $\dot{B}_{\mathrm{L}1}, \dot{B}_{\mathrm{L}2}$, which depend on the in-plane coordinates $X_1, X_2$ and are independent of $X_3$. These must satisfy (38), which reduces to

$$\dot{B}_{\mathrm{L}1,1} + \dot{B}_{\mathrm{L}2,2} = 0, \tag{70}$$

and can be satisfied by introducing the function $\psi = \psi(X_1, X_2)$ such that

$$\dot{B}_{\mathrm{L}1} = \psi_{,2}, \quad \dot{B}_{\mathrm{L}2} = -\psi_{,1}. \tag{71}$$

In the surrounding space the increment of the magnetic induction $\dot{\mathbf{B}}^\star_{\mathrm{L}}$ satisfies $(54)_1$ and hence there exists a function $\psi^\star = \psi^\star(X_1, X_2)$ such that

$$\dot{B}^\star_{\mathrm{L}1} = \psi^\star_{,2}, \quad \dot{B}^\star_{\mathrm{L}2} = -\psi^\star_{,1}. \tag{72}$$

These, combined with (51), give the Lagrangian components of the magnetic field increment $\dot{\mathbf{H}}^\star_{\mathrm{L}}$ as

$$\begin{aligned} \dot{H}^\star_{\mathrm{L}1} &= \mu_0^{-1} \left[ \left( \lambda \dot{u}^\star_{1,2} + \lambda^{-1} \dot{u}^\star_{2,1} \right) B^\star_{\mathrm{L}2} + \lambda^2 \dot{B}^\star_{\mathrm{L}1} \right], \\ \dot{H}^\star_{\mathrm{L}2} &= \mu_0^{-1} \left[ 2\lambda^{-1} \dot{u}^\star_{2,2} B^\star_{\mathrm{L}2} + \lambda^{-2} \dot{B}^\star_{\mathrm{L}2} \right], \end{aligned} \tag{73}$$

which must satisfy $(54)_2$.

The increment in the Lagrangian Maxwell stress $\dot{\mathbf{T}}^\star$ is given in (49), which results in the component equations

$$\begin{aligned} \dot{T}^\star_{11} &= \frac{\mu_0^{-1} \lambda^{-4} B_{\mathrm{L}2}}{2} \left[ B_{\mathrm{L}2} \left( \dot{u}^\star_{1,1} - 2\lambda^2 \dot{u}^\star_{2,2} \right) - 2\lambda \dot{B}^\star_{\mathrm{L}2} \right], \\ \dot{T}^\star_{21} &= \frac{\mu_0^{-1} B_{\mathrm{L}2}}{2} \left[ B_{\mathrm{L}2} \left( 2\dot{u}^\star_{1,2} + \lambda^{-2} \dot{u}^\star_{2,1} \right) + 2\lambda \dot{B}^\star_{\mathrm{L}1} \right], \\ \dot{T}^\star_{12} &= \frac{\mu_0^{-1} \lambda^{-2} B_{\mathrm{L}2}}{2} \left( B_{\mathrm{L}2} \dot{u}^\star_{1,2} + 2\lambda \dot{B}^\star_{\mathrm{L}1} \right], \\ \dot{T}^\star_{22} &= \frac{\mu_0^{-1} \lambda^{-2} B_{\mathrm{L}2}}{2} \left( 2\dot{B}^\star_{\mathrm{L}2} + \lambda \dot{u}^\star_{2,2} B_{\mathrm{L}2} \right), \\ \dot{T}^\star_{33} &= -\mu_0^{-1} \lambda^{-2} B_{\mathrm{L}2} \left( \dot{B}^\star_{\mathrm{L}2} + \lambda \dot{u}^\star_{2,2} B_{\mathrm{L}2} \right). \end{aligned} \tag{74}$$

These must satisfy the equilibrium equation (48), which reduces to the components

$$\dot{T}^\star_{11,1} + \dot{T}^\star_{21,2} = 0, \quad \dot{T}^\star_{12,1} + \dot{T}^\star_{22,2} = 0, \tag{75}$$

resulting in

$$\psi^\star_{,11} + \lambda^4 \psi^\star_{,22} = 0, \quad \psi^\star_{,21} - \psi^\star_{,12} = 0, \tag{76}$$

respectively.

In the material we consider the energy $\Omega$ as a function of the invariants (61) and (63). A calculations using (46) and $(40)_{1,2}$ results in the nonzero second and third order tensor components

$$\mathcal{K}_{11}, \mathcal{K}_{22} \quad \text{and} \quad \Gamma_{112}, \Gamma_{222}, \Gamma_{211}, \Gamma_{121}. \tag{77}$$

In addition, we have the nonzero fourth-order terms

$$\mathcal{A}_{1111}, \mathcal{A}_{2222}, \mathcal{A}_{1122}, \mathcal{A}_{2211}, \mathcal{A}_{1212}, \mathcal{A}_{2121}, \mathcal{A}_{2112}, \mathcal{A}_{1221}, \tag{78}$$

where we point to the symmetry $\mathcal{A}_{\alpha i \beta j} = \mathcal{A}_{\beta j \alpha i}$.

We now use (68) and the magnetic induction increments $\dot{B}_{L1}, \dot{B}_{L2}$ in (44) and obtain the increments of the magnetic field

$$\dot{H}_{\mathrm{L}1} = \Gamma_{211} \dot{u}_{1,2} + \Gamma_{121} \dot{u}_{2,1} + \mathcal{K}_{11} \dot{B}_{L1},$$

$$\dot{H}_{\mathrm{L}2} = \Gamma_{112} \dot{u}_{1,1} + \Gamma_{222} \dot{u}_{2,2} + \mathcal{K}_{22} \dot{B}_{L2}, \tag{79}$$

and by using (47) we obtain the connection

$$\dot{H}_{\mathrm{L}1,2} - \dot{H}_{\mathrm{L}2,1} = 0. \tag{80}$$

From (39), the nonzero components of the incremental stress $\dot{\mathbf{T}}$ are found to be

$$\begin{aligned} \dot{T}_{11} &= \left( \mathcal{A}_{1111} + p\lambda^{-2} \right) \dot{u}_{1,1} + \mathcal{A}_{1122} \dot{u}_{2,2} - \dot{p}\lambda^{-1} + \Gamma_{112} \dot{B}_{L2}, \\ \dot{T}_{22} &= \left( \mathcal{A}_{2222} + p\lambda^{2} \right) \dot{u}_{2,2} + \mathcal{A}_{2211} \dot{u}_{1,1} - \dot{p}\lambda + \Gamma_{222} \dot{B}_{L2}, \\ \dot{T}_{12} &= \left( \mathcal{A}_{1221} + p \right) \dot{u}_{1,2} + \mathcal{A}_{1212} \dot{u}_{2,1} + \Gamma_{121} \dot{B}_{L1}, \\ \dot{T}_{21} &= \left( \mathcal{A}_{2112} + p \right) \dot{u}_{2,1} + \mathcal{A}_{2121} \dot{u}_{1,2} + \Gamma_{211} \dot{B}_{L1}. \end{aligned} \tag{81}$$

These must satisfy the equilibrium condition (41), which specialized to the current condition results in

$$\dot{T}_{11,1} + \dot{T}_{21,2} = 0, \quad \dot{T}_{12,1} + \dot{T}_{22,2} = 0. \tag{82}$$

### 6.1. Boundary conditions at $X_2 = 0$

Consider an interface that divides the material in upper and lower parts distinguished as the upper layer and the substrate, respectively. The corresponding field vectors on the two sides are then identified with subscripts u and s. In addition, consider a Lagrangian unit vector $\mathbf{N}$ normal to the interface pointing from the side s to the side u. The jump in a vector on the interface is the difference between its values on side u and side s. For example, the magnetic induction $\mathbf{B}_{\mathrm{L}}$ has jump $\mathbf{B}_{\mathrm{Lu}} - \mathbf{B}_{\mathrm{Ls}}$.

Then, the incremental forms of the boundary conditions at the material interface $X_2 = 0$ are obtained from (55), (56) and (57) by replacing the contributions from the Maxwell stress $\boldsymbol{\tau}^\star$, the magnetic fields $\mathbf{B}^\star$ and $\mathbf{H}^\star$ exterior to the material by the increments of the nominal stress $\dot{\mathbf{T}}$, the magnetic induction and the magnetic field existing in the upper part of the material interface. The jump conditions of the total nominal stress $\dot{\mathbf{T}}$, of the magnetic induction $\dot{\mathbf{B}}_{\mathrm{L}}$ and of the magnetic field $\dot{\mathbf{H}}_{\mathrm{L}}$ are summarized as

$$\left( \dot{\mathbf{T}}^{\mathrm{T}}_{\mathrm{u}} - \dot{\mathbf{T}}^{\mathrm{T}}_{\mathrm{s}} \right) \mathbf{N} = \mathbf{0}, \quad \left( \dot{\mathbf{B}}_{\mathrm{Lu}} - \dot{\mathbf{B}}_{\mathrm{Ls}} \right) \cdot \mathbf{N} = 0, \quad \left( \dot{\mathbf{H}}_{\mathrm{Lu}} - \dot{\mathbf{H}}_{\mathrm{Ls}} \right) \times \mathbf{N} = \mathbf{0}, \tag{83}$$

where we assume no increments of the free surface current and of the mechanical traction $\dot{\mathbf{t}}_{\mathrm{A}} = \mathbf{0}$. In addition to the boundary conditions (83) we have the incremental displacement continuity $\dot{\mathbf{u}}_{\mathrm{s}} = \dot{\mathbf{u}}_{\mathrm{u}}$.

### 6.2. Boundary conditions at $X_2 = H$

We now specify the incremental boundary conditions (55), (56), (57) at the surface $X_2 = H$, where we take the increment of an applied mechanical traction $\dot{\mathbf{t}}_{\mathrm{A}} = \mathbf{0}$. Then, the incremental traction boundary condition (55) has the component equations

$$\dot{T}_{21} + \dot{u}^\star_{2,1} \tau^\star_{11} - \lambda \dot{\tau}^\star_{12} = 0, \quad \dot{T}_{22} + \lambda^2 \dot{u}^\star_{2,2} \tau^\star_{22} - \lambda \dot{\tau}^\star_{22} = 0. \tag{84}$$

The increments in the magnetic field and in the magnetic induction satisfy, respectively,

$$\dot{H}_{\mathrm{L}1} - \lambda \dot{H}^\star_1 + H^\star_2 \dot{u}^\star_{2,1} = 0, \quad \dot{B}_{\mathrm{L}2} - \lambda \dot{B}^\star_2 + \lambda^2 \dot{u}^\star_{2,2} B^\star_2 = 0, \tag{85}$$

respectively. As for the material interface we have the displacement continuity $\dot{\mathbf{u}}_{\mathrm{u}} = \dot{\mathbf{u}}^\star$.

## 7. Solution procedure

We seek incremental displacements that are periodic in the $X_1$ direction, depend exponentially on $X_2$ and are subject to the incom-pressibility condition (69). For the solid we take the forms

$$\dot{u}_1 = F \sin\left(KX_1\right) e^{rKX_2}, \quad \dot{u}_2 = G \cos\left(KX_1\right) e^{rKX_2}, \tag{86}$$

where $K$ is the wavenumber of the perturbation, $r$ is a parameter to be determined and $F, G$ are constants. Outside the material we take the incremental displacement components as

$$\dot{u}_1^\star = F^\star \sin\left(KX_1\right) e^{rKX_2}, \quad \dot{u}_2^\star = G^\star \cos\left(KX_1\right) e^{rKX_2}, \tag{87}$$

and impose continuity at the interfaces $X_2 = 0$ and $X_2 = H$.

To evaluate the stress components (81), in addition to the incremental displacement (86) we also need the increment of the Lagrange multiplier $\dot{p}$, which we take as

$$\dot{p} = K\lambda^{-1} P \cos\left(KX_1\right) e^{rKX_2}, \tag{88}$$

where $P$ is a constant. An alternative method is to eliminate $\dot{p}$ from the two equilibrium equations (82) by cross-differentiation thereby reducing the number of governing equations, see Dorfmann and Ogden (2010), for example.

For the increments of the magnetic induction in the material (71) and in the surrounding space (72) we use

$$\psi = V \sin\left(KX_1\right) e^{rKX_2}, \quad \psi^\star = V^\star \sin\left(KX_1\right) e^{rKX_2}, \tag{89}$$

respectively.

We are now able to expand (80) and (82) in terms of the small-amplitude solutions (86), (88) and $(89)_1$. Combined with the incompressibility condition (69) this yields four homogeneous linear equations. For a non-trivial solution to exist, the determinant of the coefficient matrix of $F, G, V, P$ must vanish resulting in the bicubic equation in $r^2$ given by

$$\begin{aligned}
&\left(\mathcal{A}_{2121}\mathcal{K}_{11} - \Gamma_{211}^2\right)\lambda^4 r^6 - \left[\mathcal{A}_{2222}\mathcal{K}_{11} - 2\left(\mathcal{A}_{2211}\mathcal{K}_{11} + \mathcal{A}_{2112}\mathcal{K}_{11} + \Gamma_{211}\Gamma_{222}\right.\right.\\
&\left.\left. -\Gamma_{121}\Gamma_{211}\right)\lambda^2 + \left(\mathcal{A}_{1111}\mathcal{K}_{11} + \mathcal{A}_{2121}\mathcal{K}_{22} + 2\Gamma_{112}\Gamma_{211}\right)\lambda^4\right] r^4 + \left[\mathcal{A}_{1212}\mathcal{K}_{11}\right.\\
&+ \mathcal{A}_{2222}\mathcal{K}_{22} - \Gamma_{121}^2 + 2\Gamma_{121}\Gamma_{222} - \Gamma_{222}^2 - 2\left(\mathcal{A}_{2112}\mathcal{K}_{22} + \mathcal{A}_{2211}\mathcal{K}_{22}\right.\\
&\left.\left. +\Gamma_{112}\Gamma_{121} - \Gamma_{112}\Gamma_{222}\right)\lambda^2 + \left(\mathcal{A}_{1111}\mathcal{K}_{22} - \Gamma_{112}^2\right)\lambda^4\right] r^2 - \mathcal{A}_{1212}\mathcal{K}_{22} = 0.
\end{aligned} \tag{90}$$

For a given energy function, Eq. (90) yields three positive and three negative solutions. In the substrate we want the perturbation to vanish for $X_2 \to -\infty$, hence from (86) follows that only positive values are retained. The general solution in the substrate is then obtained using

$$\dot{u}_{1\mathrm{s}} = \sum_{i=1}^{3} F_{\mathrm{s}i} \sin\left(KX_1\right) e^{r_i KX_2}, \quad \dot{u}_{2\mathrm{s}} = \sum_{i=1}^{3} G_{\mathrm{s}i} \cos\left(KX_1\right) e^{r_i KX_2}, \tag{91}$$

$$\dot{p}_{\mathrm{s}} = K\lambda^{-1} \sum_{i=1}^{3} P_{\mathrm{s}i} \cos\left(KX_1\right) e^{r_i KX_2}, \quad \psi_{\mathrm{s}} = \sum_{i=1}^{3} V_{\mathrm{s}i} \sin\left(KX_1\right) e^{r_i KX_2}, \tag{92}$$

where $F_{\mathrm{s}i}, G_{\mathrm{s}i}, V_{\mathrm{s}i}, P_{\mathrm{s}i}$, $i \in \{1,2,3\}$ are constants and $r_i > 0$. Eqs. (91) and (92) with $i \in \{1,2,3\}$ are now used in (69), (80) and (82) to express $G_{\mathrm{s}i}$, $P_{\mathrm{s}i}$ and $V_{\mathrm{s}i}$ in terms of $F_{\mathrm{s}i}$. This allows to write the twelve unknowns in terms of three independent quantities taken as $F_{\mathrm{s}1}, F_{\mathrm{s}2}, F_{\mathrm{s}3}$.

In the upper layer we consider all six solutions $r_i, i \in \{1,\dots,6\}$ resulting in

$$\dot{u}_{1\mathrm{u}} = \sum_{i=1}^{6} F_{\mathrm{u}i} \sin\left(KX_1\right) e^{r_i KX_2}, \quad \dot{u}_{2\mathrm{u}} = \sum_{i=1}^{6} G_{\mathrm{u}i} \cos\left(KX_1\right) e^{r_i KX_2}, \tag{93}$$

$$\dot{p}_{\mathrm{u}} = K\lambda^{-1} \sum_{i=1}^{6} P_{\mathrm{u}i} \cos\left(KX_1\right) e^{r_i KX_2}, \quad \psi_{\mathrm{u}} = \sum_{i=1}^{6} V_{\mathrm{u}i} \sin\left(KX_1\right) e^{r_i KX_2}. \tag{94}$$

We repeat the process to reduce the 24 constants $F_{\mathrm{u}i}, G_{\mathrm{u}i}, P_{\mathrm{u}i}, V_{\mathrm{u}i}$, $i \in \{1,\dots,6\}$ to a set of six, which are selected as $F_{\mathrm{u}i}$, $i \in \{1,\dots,6\}$.

Outside the material the incremental Maxwell stress $\dot{\mathbf{T}}^\star$ must satisfy the equilibrium condition (48), which has the component equations (75) and yields the constraint (76). In addition, the magnetic field components (73) must satisfy $(54)_2$. From (76) and $(89)_2$ we obtain the two roots $r = \pm\lambda^{-2}$. For the perturbation to decay for increasing values of $X_2$ we restrict attention to the negative roots and take the general solution in the form

$$\dot{u}_1^\star = F^\star \sin\left(KX_1\right) e^{-\lambda^{-2}KX_2}, \quad \dot{u}_2^\star = G^\star \cos\left(KX_1\right) e^{-\lambda^{-2}KX_2} \tag{95}$$

and

$$\psi^\star = V^\star \sin\left(KX_1\right) e^{-\lambda^{-2}KX_2}. \tag{96}$$

Writing the jump conditions at the material and material-vacuum interfaces, given in (83), (84) and (85), results in 12 homogeneous linear equations for the twelve unknowns $F_{\mathrm{s}i}, i \in \{1,2,3\}$, $F_{\mathrm{u}i}, i \in \{1,\dots,6\}$ and $F^\star, G^\star, V^\star$. The vanishing of the determinant of coefficients gives the bifurcation criterion, which is solved numerically.

### 7.1. A magnetoelastic Mooney–Rivlin material

To illustrate the results, we focus on an isotropic magnetoelastic Mooney–Rivlin material with the total energy $\Omega$ given by

$$\Omega = \frac{1}{4}\mu\left[(1+\gamma)\left(I_1 - 3\right) + (1-\gamma)\left(I_2 - 3\right)\right] + \frac{1}{2\mu_0}\left(\alpha I_4 + \beta I_5\right), \tag{97}$$

where $\mu$ is the shear modulus in the reference configuration in the absence of a magnetic field, $\gamma$ is a dimensionless material constant and $\alpha, \beta$ are magnetoelastic coupling parameters. To obtain numerical results we must differentiate between material parameters of the substrate and of the upper layer, which will be noted as necessary. Using (97) in (66) we obtain the nominal stress components as

$$\begin{aligned}
T_{11} &= \frac{1}{2}\mu\left(2\lambda + \lambda^{-1} - \gamma\lambda^{-1}\right) - p\lambda^{-1},\\
T_{22} &= \frac{1}{2}\mu\left(\lambda + 2\lambda^{-1} - \gamma\lambda\right) - p\lambda + \mu_0^{-1}\beta\lambda^{-1} I_4,\\
T_{33} &= \frac{1}{2}\mu\left[1 + \gamma + (1-\gamma)\left(\lambda^2 + \lambda^{-2}\right)\right] - p,
\end{aligned} \tag{98}$$

and in (67) the nonzero Lagrangian component of the magnetic field

$$H_{\mathrm{L}2} = \mu_0^{-1}\left(\alpha + \lambda^{-2}\beta\right) B_{\mathrm{L}2}. \tag{99}$$

The traction continuity across the interface of the layered half-space and the continuity of the component $T_{22}$ and the Maxwell stress $T_{22}^\star$ on the boundary surface $X_2 = H$ give the Lagrangian multipliers in the substrate and in the upper layer, respectively, as

$$\begin{aligned}
p_{\mathrm{s}} &= \frac{1}{2}\mu_{\mathrm{s}}\left(1 - \gamma_{\mathrm{s}} + 2\lambda^{-2}\right) + \frac{1}{2}\mu_0^{-1} I_4\left(2\beta_{\mathrm{s}} - 1\right)\lambda^{-2},\\
p_{\mathrm{u}} &= \frac{1}{2}\mu_{\mathrm{u}}\left(1 - \gamma_{\mathrm{u}} + 2\lambda^{-2}\right) + \frac{1}{2}\mu_0^{-1} I_4\left(2\beta_{\mathrm{u}} - 1\right)\lambda^{-2}.
\end{aligned} \tag{100}$$

We now evaluate the tensor components (77) and (78) to specialize (90) to

$$\left(r^2 - 1\right)\left(r^2\lambda^4 - 1\right)\left[r^2\left(\alpha + \alpha\beta\bar{I}_4 + \beta\lambda^2\right)\lambda^2 - \alpha\lambda^2 - \beta\right] = 0, \tag{101}$$

where we introduced the dimensionless measure $\bar{I}_4$ defined by

$$\bar{I}_4 = \bar{B}_{\mathrm{L}2}^2, \quad \bar{B}_{\mathrm{L}2} = \frac{B_{\mathrm{L}2}}{\sqrt{\mu\mu_0}}. \tag{102}$$

From (101) we find the roots

$$r_{1,4} = \pm 1, \quad r_{2,5} = \pm\lambda^{-2}, \quad r_{3,6} = \pm\frac{\sqrt{\alpha\lambda^2 + \beta}}{\sqrt{\left(\alpha + \alpha\beta\bar{I}_4 + \beta\lambda^2\right)\lambda^2}}, \tag{103}$$

which are used to construct the general solution in the substrate (91), (92) and in the upper layer (93) and (94), where the specified material parameters $\mu, \alpha, \beta$ will be denoted $\mu_{\mathrm{s}}, \alpha_{\mathrm{s}}, \beta_{\mathrm{s}}$ and $\mu_{\mathrm{u}}, \alpha_{\mathrm{u}}, \beta_{\mathrm{u}}$, respectively.

To write the governing equations and boundary conditions in dimensionless form we define the constants

$$\bar{F} = \frac{F}{H}, \quad \bar{G} = \frac{G}{H}, \quad \bar{V} = \frac{V}{H\sqrt{\mu\mu_0}}, \quad \bar{P} = \frac{P}{\mu}, \quad \bar{H} = 1, \tag{104}$$

where $H$ is the undeformed thickness of the upper layer. Using (102) in (100) gives the dimensionless forms of the Lagrangian multipliers.

A similar process is used to obtain the dimensionless quantities $\bar{F}^{\star}, \bar{G}^{\star}$ and $\bar{V}^{\star}$.

The constants in the substrate $\bar{F}_{si}, \bar{G}_{si}, \bar{V}_{si}, \bar{P}_{si}, i \in \{1,2,3\}$ are not independent and are connected via the incompressibility condition (69) given by

$$\bar{F}_{si} + r_{si}\lambda^2 \bar{G}_{si} = 0, \tag{105}$$

where the subscript $i \in \{1,2,3\}$ and the summation over repeated indexes does not apply. They are also connected via (80)

$$r_{si}^2 \beta_s \lambda^3 \bar{B}_{12} \bar{F}_{si} + r_{si}\beta_s \lambda \bar{B}_{L2} \bar{G}_{si} + \left[\alpha_s \left(r_{si}^2 - 1\right)\lambda^2 + \beta_s \left(r_{si}^2 \lambda^4 - 1\right)\right] \bar{V}_{si} = 0, \tag{106}$$

and by the equilibrium equations (82), which specialize to

$$r_{si}\left[2(\gamma - 1)\lambda^2 + \left(1 - 2\beta_s\right)\bar{I}_4 - 2\right]\bar{F}_{si} + \left[2\left(1 - 2r_{si}^2\right)\lambda^2 + \left(1 - 4\beta_s\right) r_{si}^2 \lambda^2 \bar{I}_4 + 2(\gamma - 1) r_{si}^2 \lambda^4\right]\bar{G}_{si} + 2 r_{si}\beta_s \lambda \bar{V}_{si} + 2 r_{si}\lambda^2 \bar{P}_{si} = 0$$

and

$$\left\{2 + 2(1-\gamma)\lambda^2 + 2\left(1 - r_{si}^2\right)\lambda^4 - \left[1 + 2\beta_s\left(\lambda^4 r_{si}^2 - 1\right)\right]\bar{I}_4\right\}\bar{F}_{si} + r_{si}\lambda^2\left[2 + \left(2\beta_s - 1\right)\bar{I}_4 - 2(\gamma - 1)\lambda^2\right]\bar{G}_{si} - 2 r_{si}^2 \beta_s \lambda^5 \bar{B}_{12}\bar{V}_{si} - 2\lambda^2 \bar{P}_{si} = 0. \tag{107}$$

These are used to write the twelve unknowns in terms of the three independent quantities $\bar{F}_{s1}, \bar{F}_{s2}, \bar{F}_{s3}$. We repeat an identical process for the upper layer and express the 24 constants $\bar{F}_{ui}, \bar{G}_{ui}, \bar{V}_{ui}, \bar{P}_{ui}, i \in \{1,\dots,6\}$ in terms of the six quantities $\bar{F}_{ui}, i \in \{1,\dots,6\}$. Adding the three constants $\bar{F}^{\star}, \bar{G}^{\star}, \bar{V}^{\star}$ valid for the outside space, the total number of dimensionless constants is reduced to twelve.

In (83), the unit vector $\mathbf{N}$ has a single nonzero component in the $X_2$ direction. Then, at $X_2 = 0$ the traction boundary condition reduces to two component equations, the continuity of the magnetic induction and the magnetic field increments simplify to one equation each. We also impose the continuity of the incremental displacement components resulting in a total of six homogeneous equations. A similar situation arises at the boundary $X_2 = H$, where two traction boundary conditions, the continuity of the normal component of the magnetic induction and of the tangent component of the magnetic field increments, combined with the displacement continuity provide the additional six equations. In summary, therefore, the jump conditions result in 12 homogeneous linear equations in the twelve unknowns $\bar{F}_{si}, i \in \{1,2,3\}$, $\bar{F}_{ui}, i \in \{1,\dots,6\}$ and $\bar{F}^{\star}, \bar{G}^{\star}, \bar{V}^{\star}$. For a known applied field and for given material parameters, we seek the value of $\lambda$ for which the determinant of coefficients vanishes.

## 8. Numerical results

To illustrate the theory, we examine four distinct scenarios as illustrated in Fig. 1. First, we verify the critical stretch, $\lambda_{cr}$, for surface instability of a layered elastic half-space composed of non-magnetizable material, with layers characterized by different stiffness ratios. Next, we analyze deformation and magnetic field increments superimposed on a homogeneously deformed magnetoelastic half-space under the influence of a magnetic field normal to its upper boundary. Third, we study the response of a non-magnetizable elastic substrate combined with a magnetizable upper layer, again subjected to deformation and magnetic field increments. Finally, we determine the bifurcation criterion for a layered magnetoelastic solid in which the substrate and the upper layer exhibit differing mechanical and magnetic properties, once more under a magnetic field normal to the upper boundary.

While there are asymptotic expressions for the limit of layered systems in the purely mechanical case (Alawiye et al., 2019; Cai and Fu, 1999; Fu and Ciarletta, 2015; Fu et al., 2021, 2024), there is currently no corresponding asymptotic solution for magnetoelastic layered systems due to the strong coupling between deformation and magnetic fields and the presence of multiple eigenbranches (physical and non-physical mathematical solutions). The arc-length continuation approach presented here reproduces previous asymptotic results in the purely mechanical limit, verifying their accuracy for both critical wave number and bifurcation point.

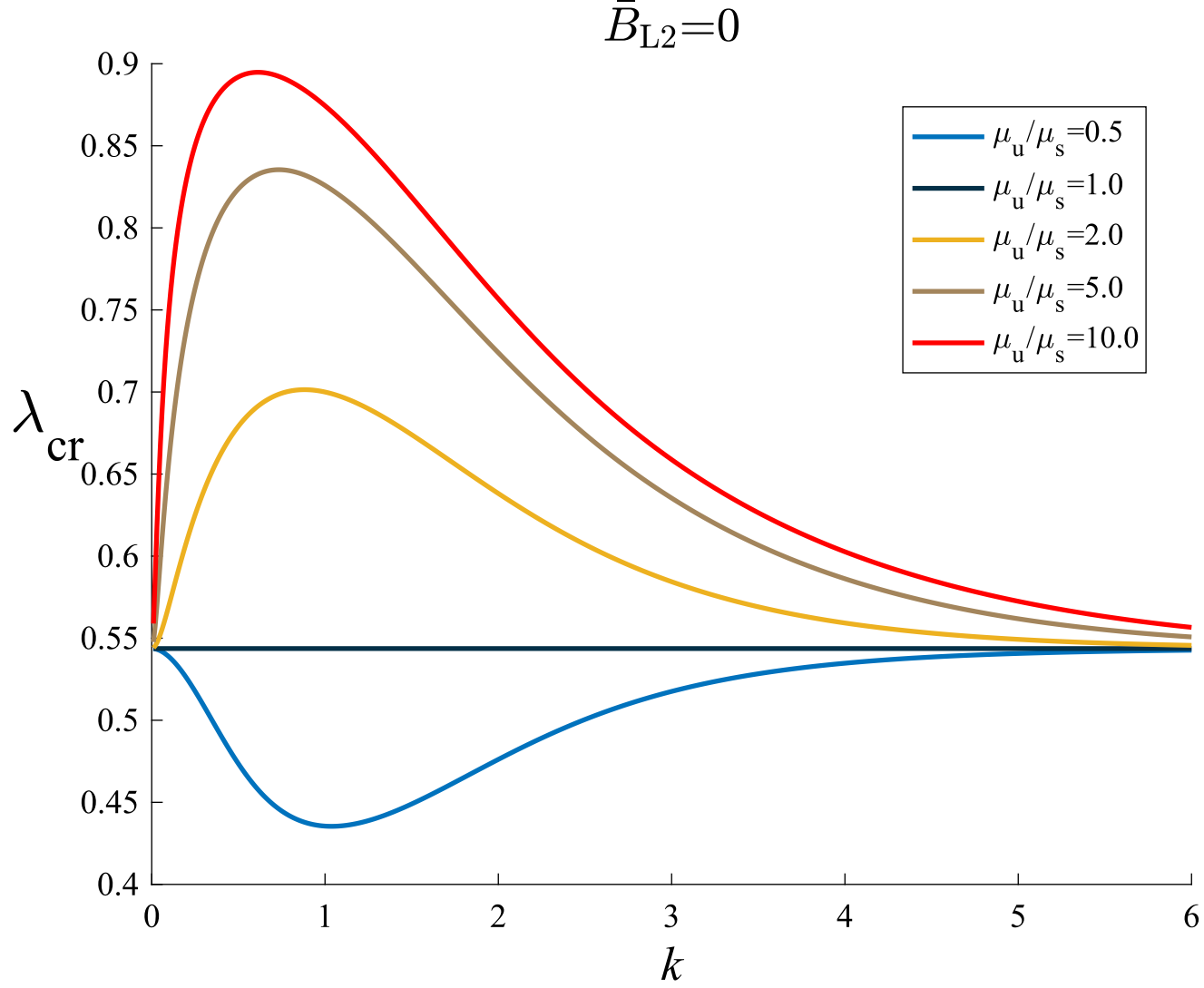

**Fig. 2.** The curves depict the critical stretch $\lambda_{cr}$ of a layered elastic half-space as a function of the wave number $k = \lambda^{-1}K$ for different stiffness ratios $\mu_u/\mu_s$. They show that surface instability of a Mooney–Rivlin material in plane strain occurs when the critical stretch becomes $\lambda_{cr} = 0.5437$, known as the Biot stretch (Biot, 1965). We find that an increase in the ratio $\mu_u/\mu_s$ destabilizes the half-space, which matches the asymptotic predictions in Alawiye et al. (2019).

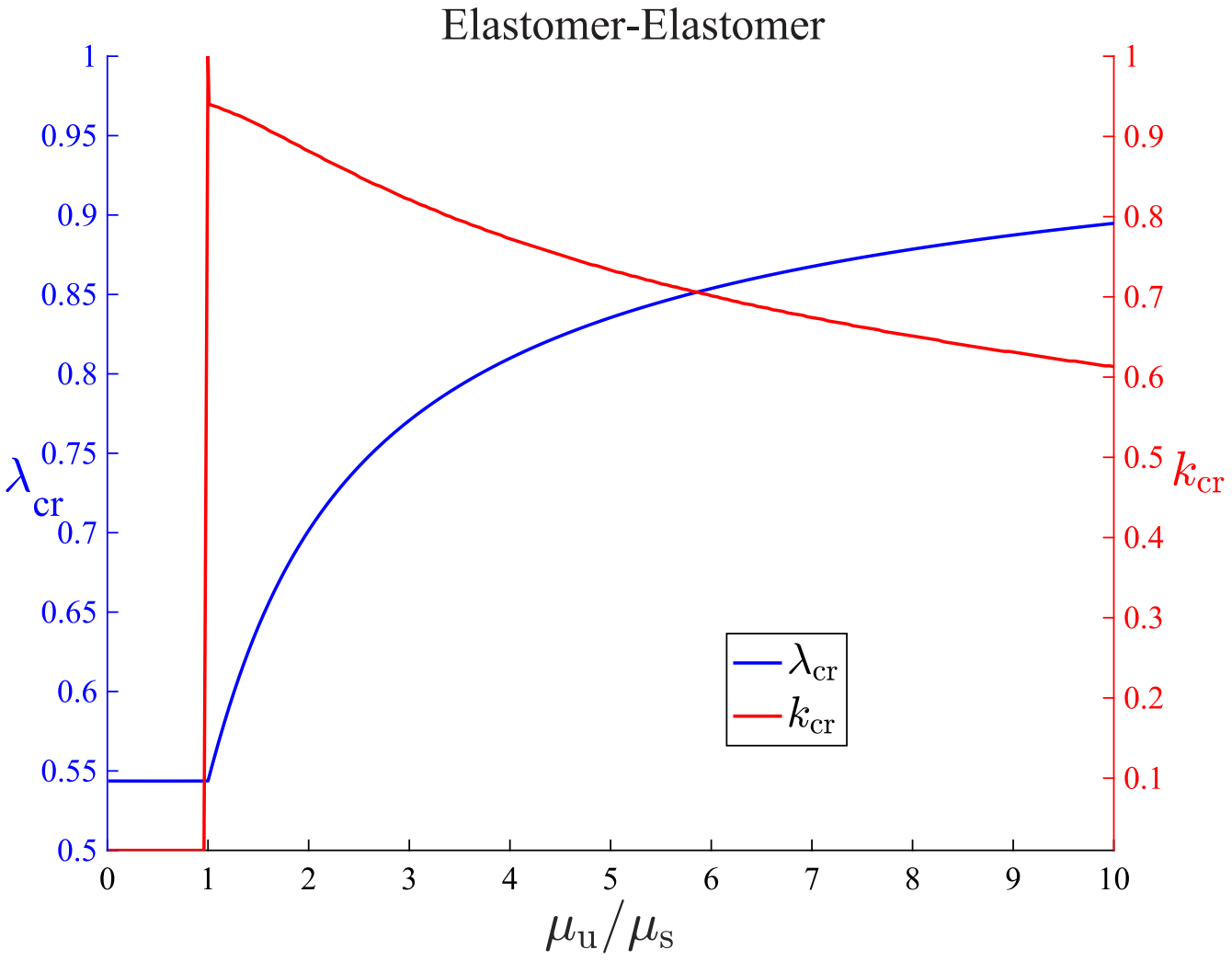

**Fig. 3.** A dual-axis plot showing on the left the critical stretch $\lambda_{cr}$ and on the right the critical wave number $k_{cr}$ as a function of the stiffness ratio $\mu_u/\mu_s$ for a layered elastic half-space with no applied magnetic field.
*Source:* Results are in agreement with prior results obtained from asymptotic predictions (Alawiye et al., 2019).

### *8.1. Layered elastic half-space of non-magnetizable material*

Fig. 2 shows the critical stretch $\lambda_{cr}$ of a layered elastic half-space as a function of the Eulerian wave number $k = \lambda^{-1}K$. The curves from top to bottom correspond to the ratios of $\mu_u/\mu_s = 0.5, 1, 2, 5, 10$. These show that for values of $k$ less than approximately 6, a stiffness ratio $\mu_u/\mu_s > 1.0$ destabilizes the half-space. This is not surprising since more of the induced load is carried by the upper layer. Also, note that for $\mu_u = \mu_s$ the critical stretch $\lambda_{cr} = 0.5437$ is independent of the wavenumber $k$. We obtain the critical stretch for surface instability of a Mooney–Rivlin material in plane strain conditions, which is known as the Biot stretch $\lambda_{cr} = 0.5437$, see Biot (1965). For a layered structure with $\mu_u/\mu_s = 10$, for example, bifurcation occurs when $\lambda_{cr} \approx 0.89$ and $k_{cr} \approx 0.61$, as shown by the red curve in Fig. 2.

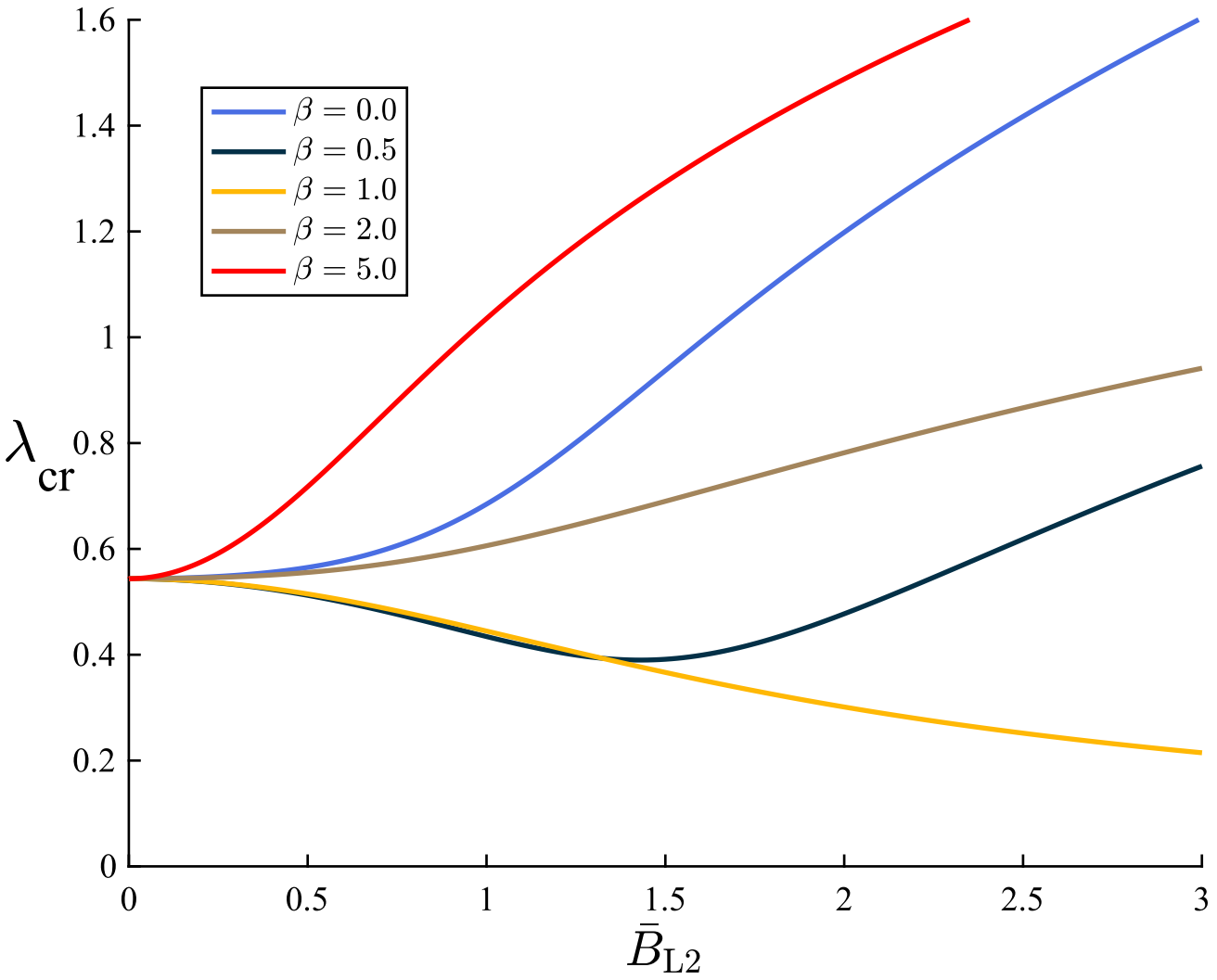

**Fig. 4.** The critical stretch $\lambda_{cr}$ of a magnetoelastic half-space as a function of $\bar{B}_{L2}$ for $\alpha = 0.5$ and $\beta = 0, 0.5, 1, 2, 5$. For increasing values of $\bar{B}_{L2}$, the half-space becomes more and more unstable, with $\beta = 1$ being an exception. In the purely mechanical case when $\bar{B}_{L2} = 0$, we recover the critical stretch of an elastic half-space $\lambda_{cr} = 0.5437$.

With the aid of the arc-length continuation method, which is capable of tracking all physical and non-physical solutions (Shahsavari and Saxena, 2025), we applied the simultaneous optimization process across different stiffness ratios presented in Fig. 2, which systematically extracted the maximum $\lambda$ value from each stability curve to identify the critical bifurcation points. The resulting critical stretch parameters $\lambda_{cr}$ and their associated critical wavenumbers $k_{cr}$ are then compiled and presented in Fig. 3, demonstrating the coupled relationship between these critical values. Fig. 3 presents a dual-axis plot showing the critical stretch $\lambda_{cr}$ and the critical wave number $k_{cr}$ as a function of the stiffness ratio $\mu_u/\mu_s$. The left axis (in blue) indicates the critical stretch, while the right axis (in red) corresponds to the critical wave number. In particular for $\mu_u/\mu_s < 1$, the critical stretch remains fixed at $\lambda_{cr} = 0.5437$ (the classical Biot solution), and the wave number reduces to zero. For $\mu_u/\mu_s = 1$, the same value $\lambda_{cr} = 0.5437$ is obtained, independent of $k$. As the stiffness ratio increases, both quantities vary, and for $\mu_u/\mu_s = 10$, for example, one finds $\lambda_{cr} \approx 0.89$ and $k_{cr} \approx 0.61$, consistent with the values shown in Fig. 2.

### 8.2. A magnetoelastic half-space

We now consider a magnetoelastic half-space in the presence of a uniform magnetic field, which prior to the superimposed increments is normal to its boundary surface. In particular, we consider the isotropic Mooney–Rivlin magnetoelastic material (97) with shear stiffness $\mu$ and with magnetoelastic coupling parameters $\alpha, \beta$. The critical compression stretch for surface instability is given in Otténio et al. (2008) for $\alpha = 0.5$ and $\beta = 0, 0.5, 1, 2, 5$. We select the same material parameters to verify the formulation and validate the numerical results. These are illustrated in Fig. 4 as the critical stretch $\lambda_{cr}$ for increasing values of the dimensionless magnetic induction $\bar{B}_{L2}$. In the purely mechanical case, when $\bar{B}_{L2} = 0$, we recover the critical stretch for an elastic half-space $\lambda_{cr} = 0.5437$. The graphs show that for $\beta = 1$ an increase in $\bar{B}_{L2}$ has a stabilizing effect. For $\beta = 0.5$ an increase in $\bar{B}_{L2}$ first stabilizes and then, for further increase, destabilizes the configuration. For $\beta = 0, 2, 5$ an increase in the magnetic field renders the half-space more unstable. The figure also shows that for a magnetic induction $\bar{B}_{L2}$ larger than some critical value the half-space becomes unstable even in tension.

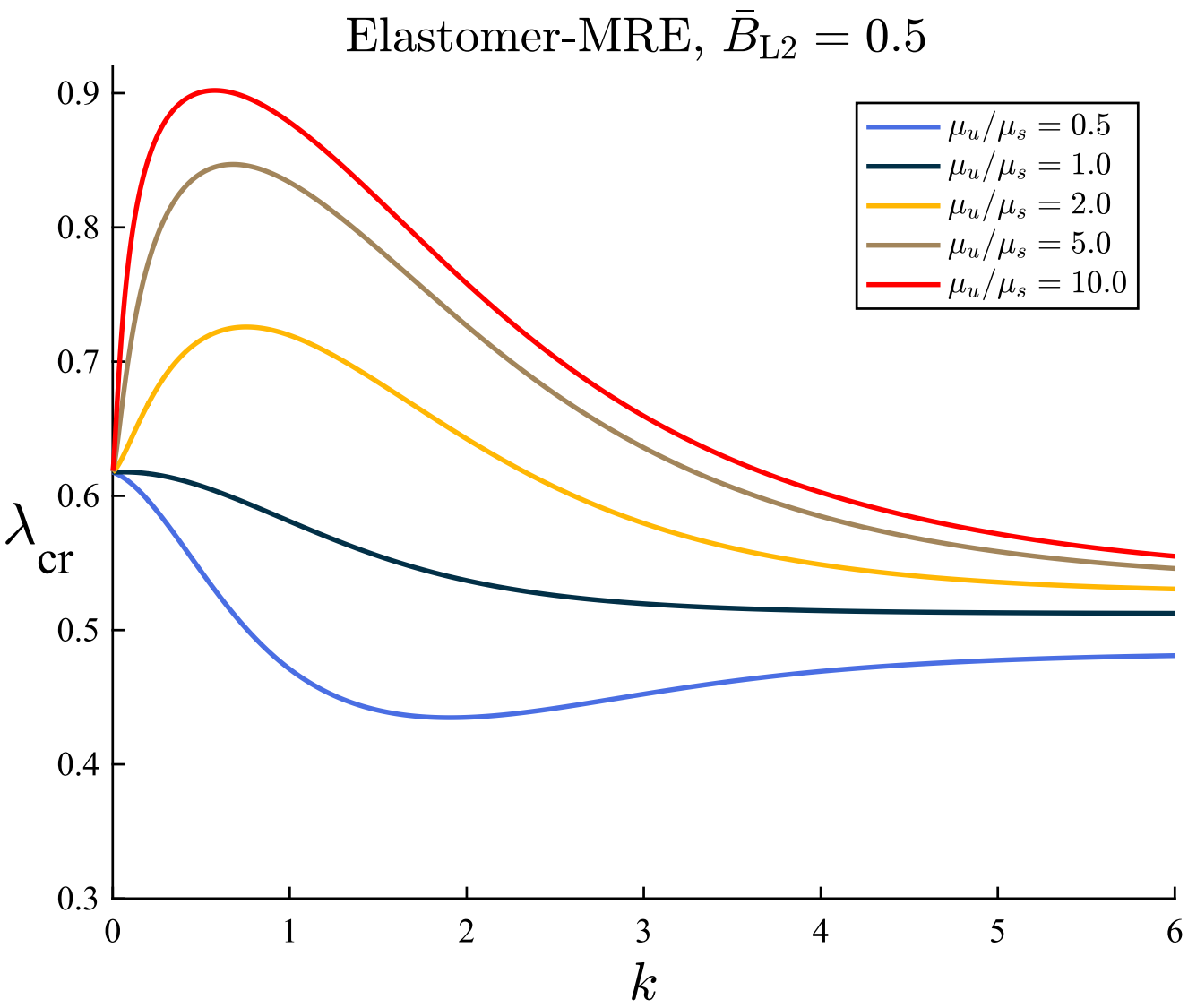

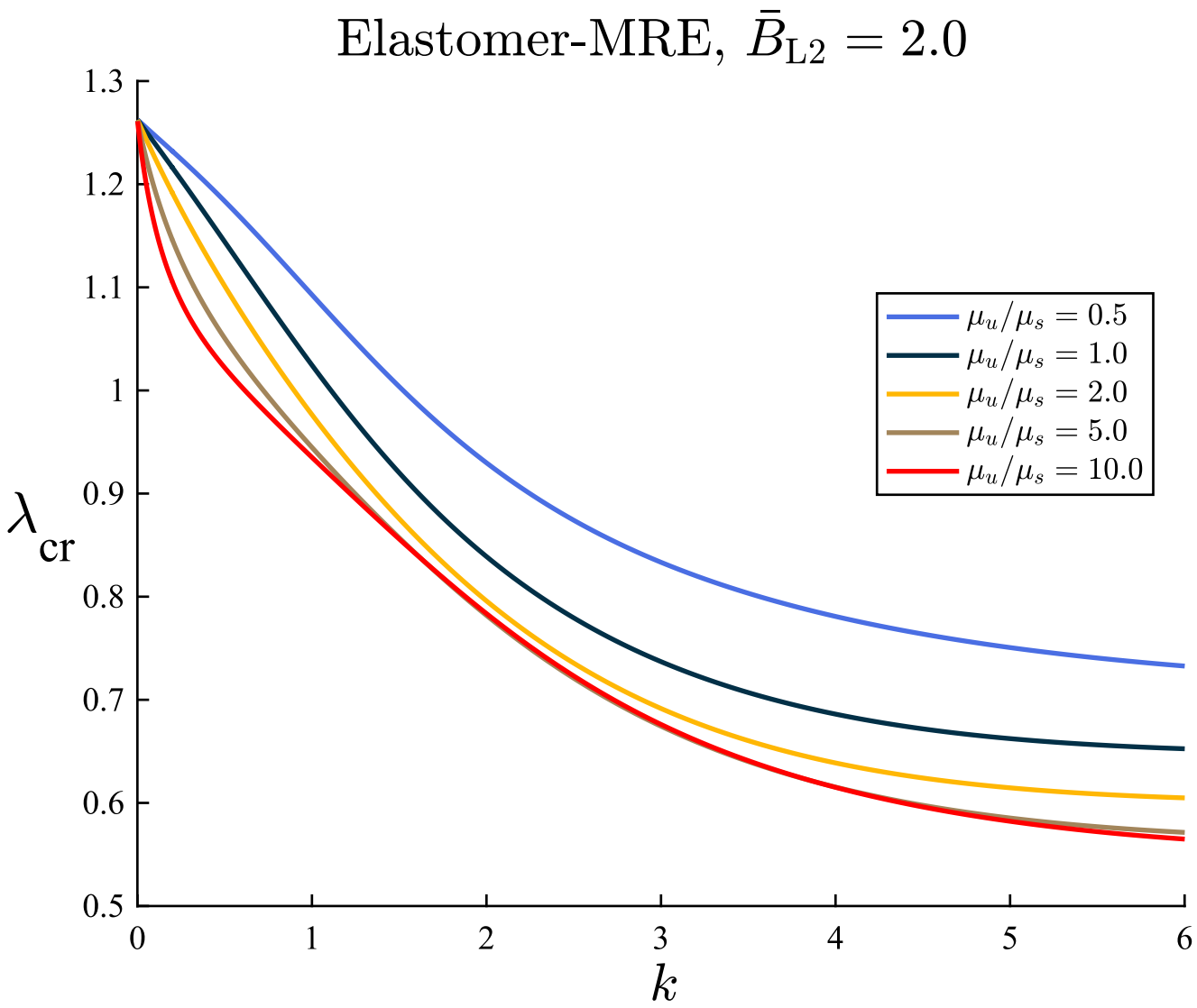

**Fig. 5.** The critical stretch $\lambda_{cr}$ of a non-magnetizable elastomer combined with a magneto-sensitive upper layer with $\alpha = \beta = 0.5$. The curves show that instability depends on the stiffness ratio and on the magnitude of $\bar{B}_{L2}$.

### 8.3. A non-magnetizable substrate with a magnetoelastic upper layer

We aim to determine the critical values $\lambda_{cr}$ and $k_{cr}$ of a layered half-space consisting of an elastic substrate and a magnetizable upper layer, as a function of the magnetic field $\bar{\mathbf{B}}_L$ applied normal to the top surface, for stiffness ratios $\mu_u/\mu_s = 0.5, 1, 2, 5, 10$. In Fig. 5, the results in the upper and lower panels correspond to $\alpha = 0.5$, $\beta = 0.5$, with $\bar{B}_{L2} = 0.5$ and $\bar{B}_{L2} = 2$, respectively. These should be compared with the response of the layered purely elastic half-space shown in Fig. 2.

The critical bifurcation values are obtained by decreasing $\lambda$ from $\lambda = 1$ until the corresponding bifurcation curve is reached. For example, for $\bar{B}_{L2} = 0.5$ and $\mu_u/\mu_s = 10$ we find $\lambda_{cr} \approx 0.90$ and $k_{cr} \approx 0.57$, the compression stretch being less than for the purely elastic half-space. Thus, less compression is required for bifurcation, indicating that the magnetic field has a destabilizing effect. Similar behavior is observed for $\mu_u/\mu_s = 0.5, 1, 2, 5$. In particular, when $\mu_u/\mu_s \leqslant 1$, bifurcation occurs when $\lambda_{cr} \approx 0.61$, compared to $\lambda_{cr} = 0.5437$ for the elastic case. Again, less compression is allowed, confirming that the magnetic field destabilizes the configuration.

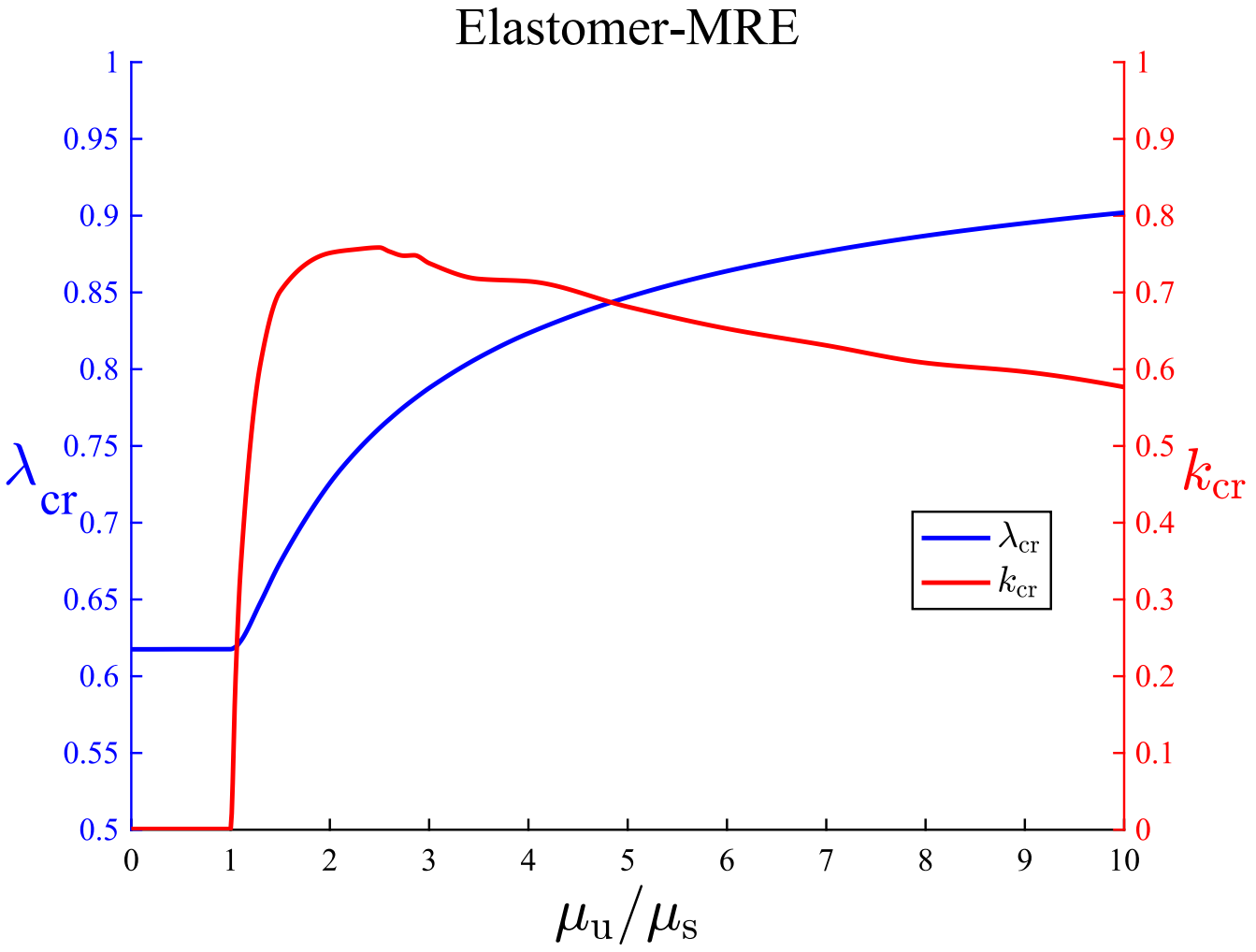

**Fig. 6.** A dual-axis plot showing on the left the critical stretch $\lambda_{cr}$ and on the right the critical wave number $k_{cr}$ as a function of the stiffness ratio $\mu_u/\mu_s$ for a non-magnetizable substrate combined with a magnetoelastic upper layer. The curves correspond to an applied field $\bar{B}_{L2} = 0.5$ and to the material constants $\alpha = \beta = 0.5$.

The response changes completely when the magnetic field $\bar{B}_{L2}$ is increased above a critical value, as illustrated in the lower panel in Fig. 5. For example, consider the undeformed configuration $\lambda = 1$ and $\mu_u/\mu_s = 10$. The results show that bifurcation occurs with a wave number $k_{cr} \approx 0.7$. For the same ratio $\mu_u/\mu_s = 10$ and for the same field $\bar{B}_{L2} = 2$, the wave number increases with the compression $\lambda$. For $\lambda = 0.8$, bifurcation occurs with $k \approx 1.9$. Comparing the results shown in the upper and lower panels we find that an increase in the magnetic field has a destabilizing effect.

In Fig. 6 we again present a dual-axis plot showing the values $\lambda_{cr}$ on the left and $k_{cr}$ on the right, as a function of the stiffness ratio $\mu_u/\mu_s$ for an applied field $\bar{B}_{L2} = 0.5$, see results reported in the upper panel in Fig. 5. Fig. 6 shows that $\mu_u/\mu_s \leqslant 1$, the stretch $\lambda_{cr} = 0.617$ and constant, while the wave number vanishes. For $\mu_u/\mu_s > 1$, the critical compression $\lambda_{cr}$ continuously reduces with the stiffness ratio reaching $\lambda_{cr} = 0.894$ for $\mu_u/\mu_s = 10$, which is identical to the results shown in Fig. 5. The wave number first increases rapidly to $k_{cr} = 0.7584$ for $\mu_u/\mu_s \approx 2.5$, the corresponding stretch being $\lambda_{cr} = 0.762$. For a further increase in the stiffness ratio, the wave number gradually decreases to $k_{cr} = 0.576$ when $\mu_u/\mu_s = 10$.

In the upper and lower panels in Fig. 5 we report the critical stretch for different stiffness ratios $\mu_u/\mu_s$ and $\bar{B}_{L2} = 0.5$ and $\bar{B}_{L2} = 2.0$, respectively. In Fig. 6 we show $\lambda_{cr}$ and $k_{cr}$ as a function of the stiffness ratio for an applied field $\bar{B}_{L2} = 0.5$, in both cases $\alpha = \beta = 0.5$. Now we illustrate the response when the applied magnetic induction assumes the values $\bar{B}_{L2} = 0, 0.5, 1, 1.5, 2, 2.5, 3$ for ratios of $\mu_u/\mu_s = 0.5, 1, 2, 5, 10$ and $\alpha = \beta = 0.5$. The corresponding critical stretch $\lambda_{cr}$ and critical wave number $k_{cr}$ are given in the upper and lower panels in Fig. 7.

To discuss the results in Fig. 7, we start with $\bar{B}_{L2} = 0$ and find that the values of $\lambda_{cr}$ and $k_{cr}$ correspond to those in Fig. 2. The values for $\bar{B}_{L2} = 0.5$ are identical to the those shown in the upper panel in Fig. 5. The findings for $\bar{B}_{L2} = 1.0$ follow the same trend. In particular, we find that with the increase of the magnetic induction the amount of compression to initiate instability is reduced, the associated wave numbers reduce as well. For $\bar{B}_{L2} \geqslant 1.5$ the behavior changes and bifurcation occurs even in tension and the critical wave number vanishes, see the insert in the upper panel of Fig. 7.

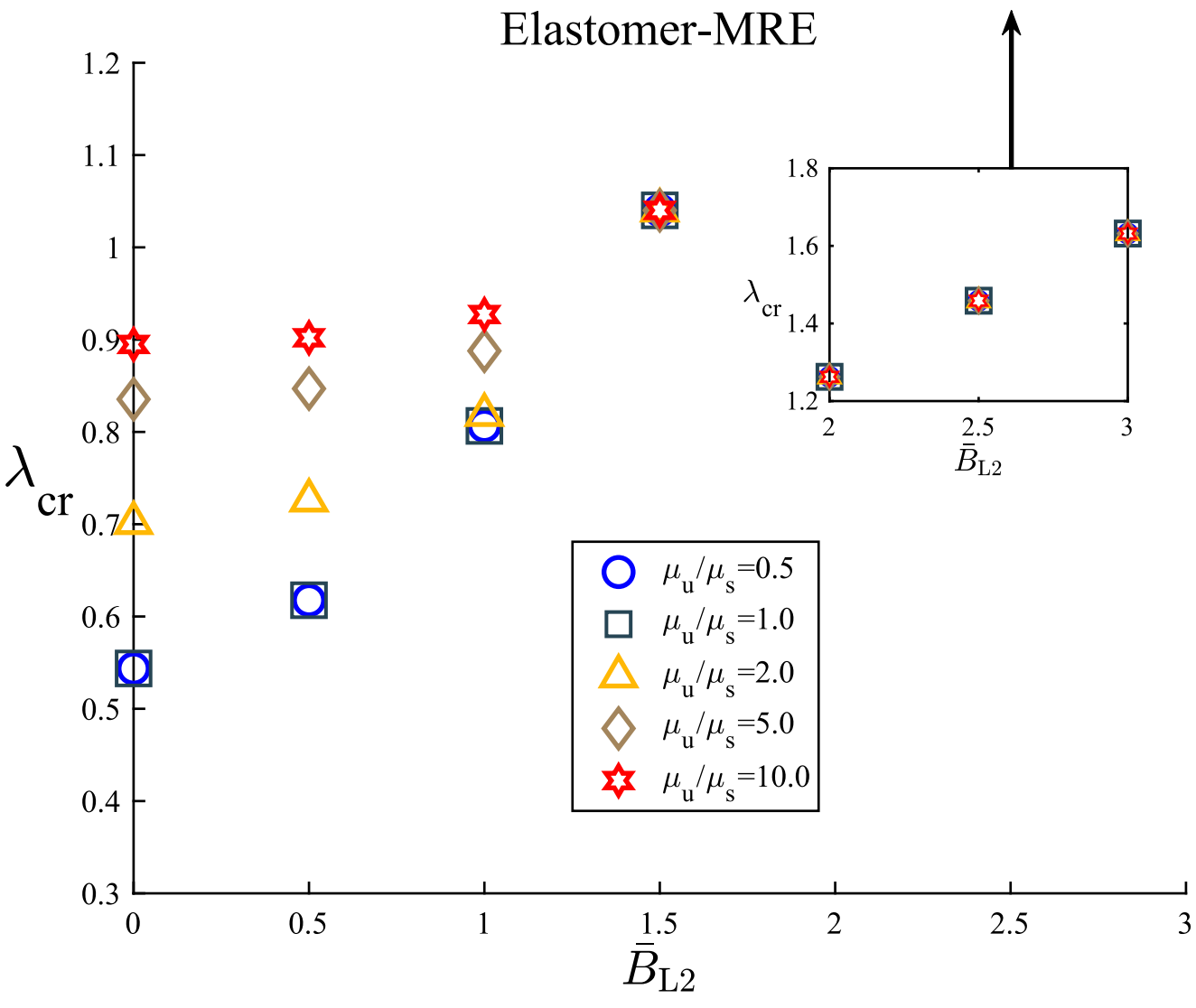

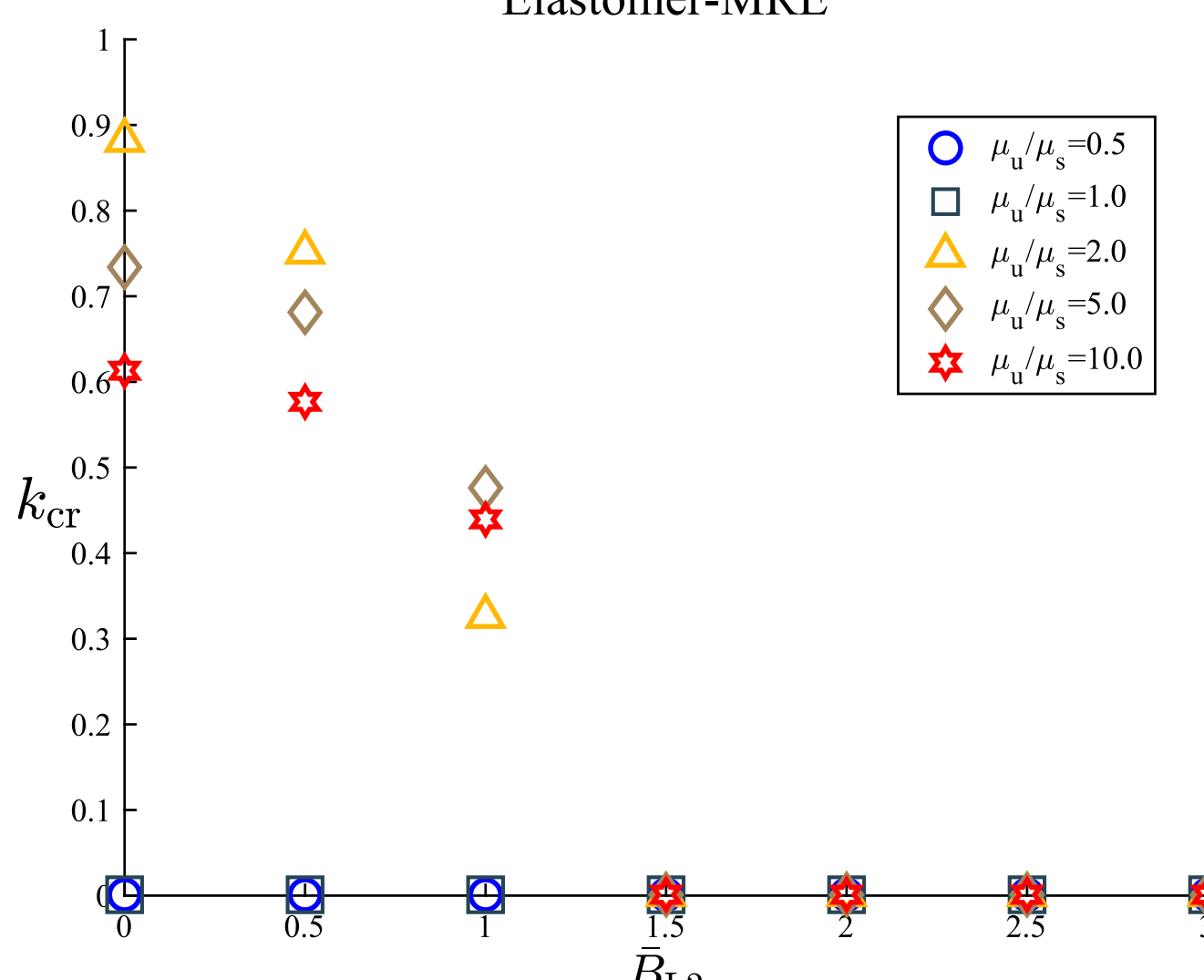

**Fig. 7.** The critical stretch $\lambda_{cr}$ (upper panel) and the critical wave number $k_{cr}$ (lower panel) of a non-magnetizable substrate combined with a magnetoelastic upper layer for $\bar{B}_{L2} = 0, 0.5, 1, 1.5, 2, 2.5, 3$, $\mu_u/\mu_s = 0.5, 1, 2, 5, 10$ and $\alpha = \beta = 0.5$.

In Fig. 5 we show the critical stretch $\lambda_{cr}$ of a non-magnetizable elastomer combined with a magneto-sensitive upper layer for $\alpha = \beta = 0.5$. Results show that surface instability depends on the stiffness ratio and on the magnitude of the applied magnetic induction. We now aim to evaluate the effect of the magnetic properties of the upper layer on the critical stretch. To accomplish this task, we consider the ratio $\mu_u/\mu_s = 1$ fixed, $\alpha = 0.5$, $\beta = 0, 0.5, 1, 2, 5$ and report the results in Fig. 8. There, the upper panel corresponds to $\bar{B}_{L2} = 0.5$, the lower to $\bar{B}_{L2} = 2$. Comparing to Fig. 5 we find that the magnetic properties can be tuned to impact the magnitude of the stretch when surface instability occurs.

### 8.4. A two-layers magnetoelastic half-space

Here we focus on the bifurcation conditions of a two-layer magnetoelastic half-space characterized by the material parameters $\alpha_s = \alpha_u = \alpha = 0.5$ and $\beta_s = \beta_u = \beta = 0.5$ and the stiffness ratios $\mu_{\mathrm{u}}/\mu_{\mathrm{s}} = 0.5, 1, 2, 5, 10$. In the upper panel of Fig. 9 we show the critical stretch and the corresponding wave number for $\bar{B}_{\mathrm{L2}} = 0.5$. These should be compared to the stability conditions in the purely elastic case and when the elastomer substrate is combined with a magnetoelastic upper layer, shown in Fig. 2, and in the upper panel of Fig. 5, respectively. In particular, for $\mu_{\mathrm{u}}/\mu_{\mathrm{s}} = 10$ we now obtain $\lambda_{\mathrm{cr}} = 0.89$ and $k_{\mathrm{cr}} = 0.61$. This indicates that a high stiffness ratio essentially eliminates the effect of the magnetic field. The impact of the magnetic coupling becomes more pronounced by reducing the stiffness ratio. For $\mu_{\mathrm{u}}/\mu_{\mathrm{s}} \leqslant 1$, for example, the two-layers magnetoelastic half-space is clearly more stable compared to the conditions shown in Figs. 2 and 5. The lower panel of Fig. 9 summarizes the stability condition when the applied field is increased to $\bar{B}_{\mathrm{L2}} = 2$. It shows that the use of a magnetoelastic substrate combined with an increase in the magnetic field stabilizes the configuration, which is true for any of the considered stiffness ratios. This also applies when compared to the conditions shown in Fig. 5.

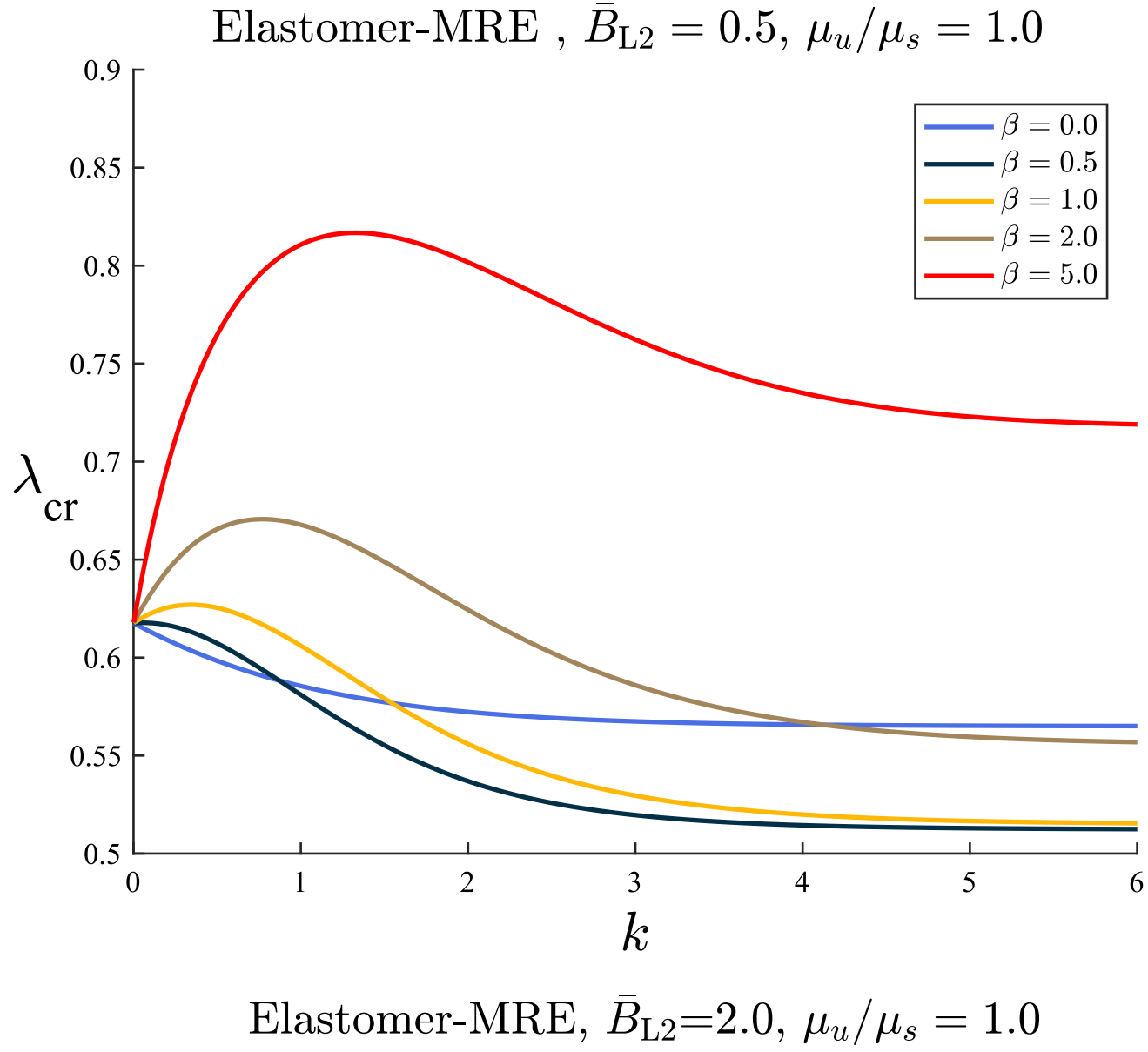

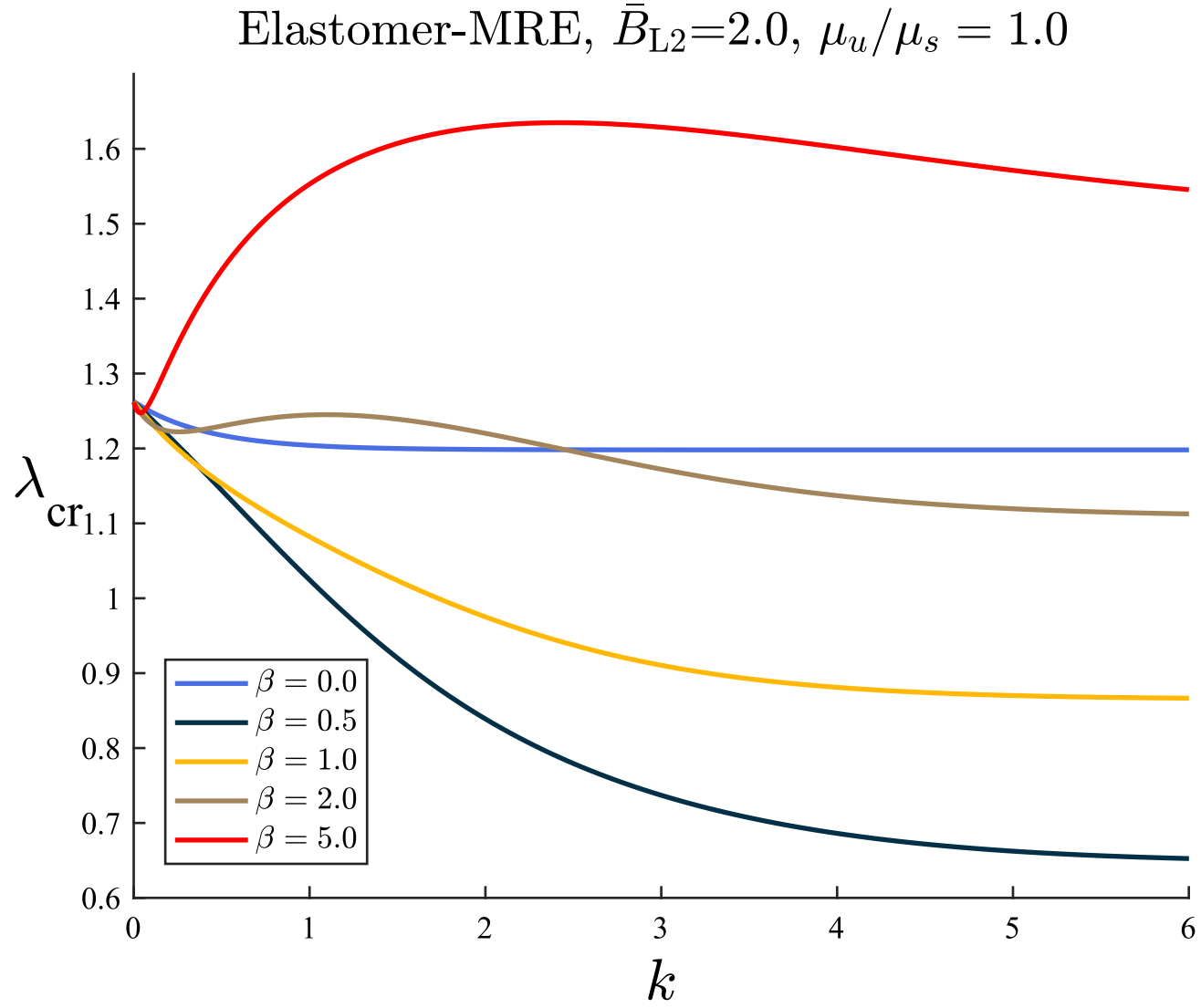

**Fig. 8.** The critical stretch $\lambda_{\mathrm{cr}}$ of a non-magnetizable elastomer combined with a magneto-sensitive upper layer with $\alpha = 0.5$ and $\beta = 0, 0.5, 1, 2, 5$. Results for $\bar{B}_{\mathrm{L2}} = 0.5$ (upper panel) and $\bar{B}_{\mathrm{L2}} = 2$ (lower panel) show that the critical stretch depends on the magnetic properties and on the magnitude of the applied magnetic field.

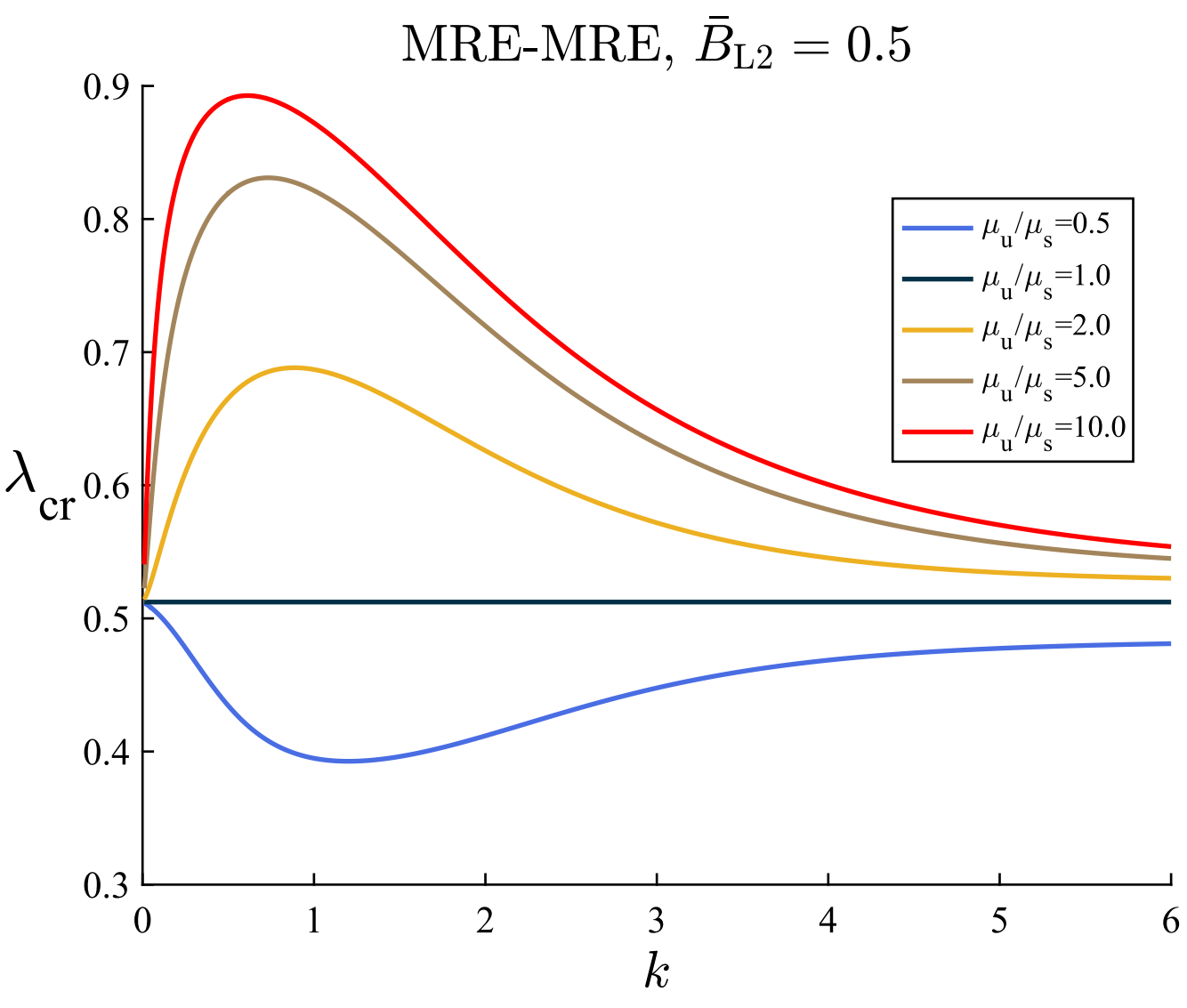

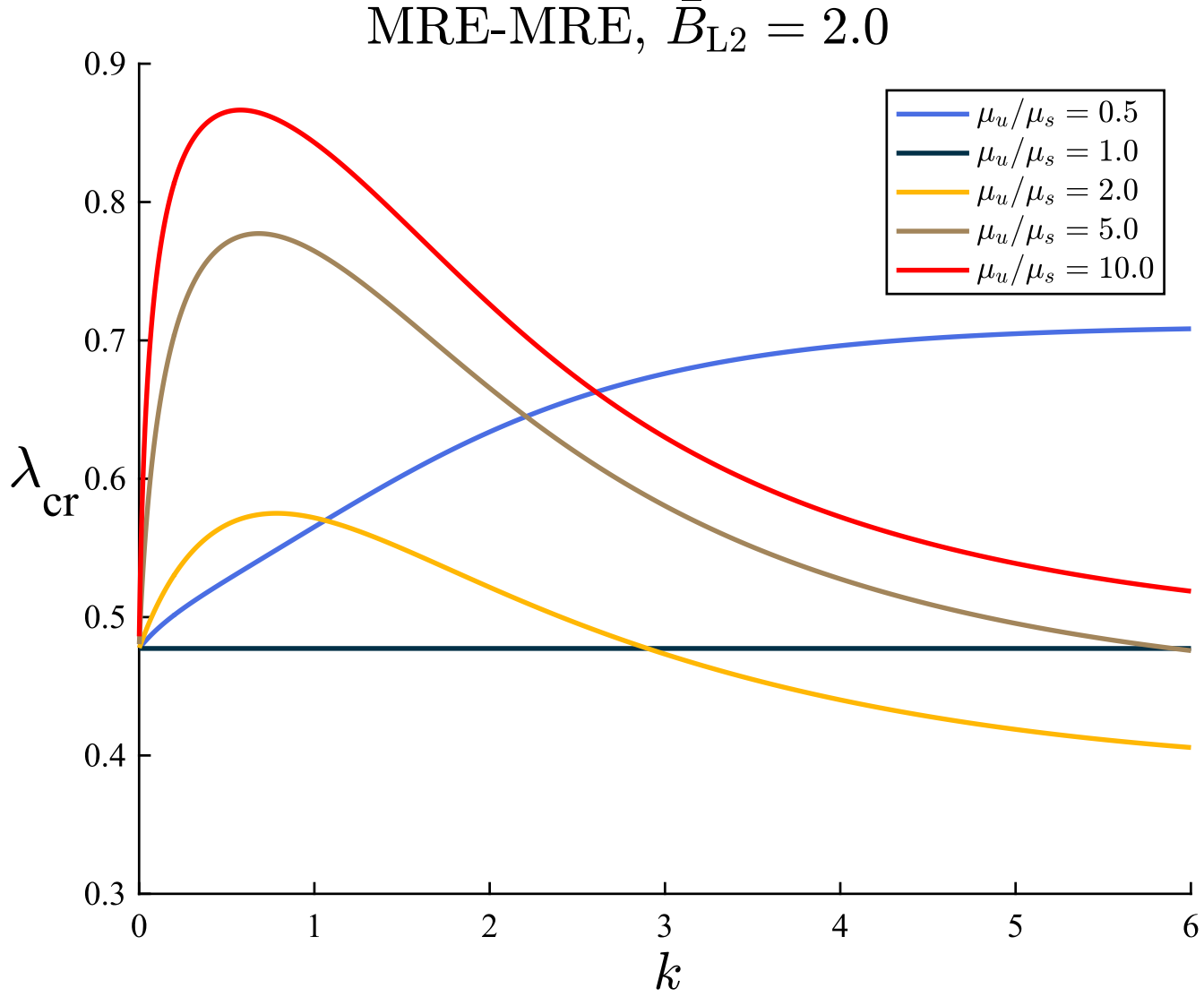

**Fig. 9.** The critical stretch $\lambda_{\mathrm{cr}}$ of a two-layers magnetoelastic half-space with $\alpha = \beta = 0.5$. The results show that the use of a magnetoelastic substrate combined with an increase in the magnetic field stabilizes the configuration, for all of the stiffness ratios considered.

The dual-axis plot in Fig. 10 illustrates the critical stretch $\lambda_{\mathrm{cr}}$ and the critical wave number $k_{\mathrm{cr}}$ when $\alpha = \beta = 0.5$ and the applied field $\bar{B}_{\mathrm{L2}} = 0.5$, as a function of $\mu_{\mathrm{u}}/\mu_{\mathrm{s}}$, see upper panel in Fig. 9. Compared to the bifurcation conditions of the layered elastic half-space in Fig. 3 and of the elastic substrate combined with a magnetoelastic upper layer in Fig. 6, we find that the magneto-sensitive substrate has a stabilizing effect for $\mu_{\mathrm{u}}/\mu_{\mathrm{s}} \leqslant 1$. For a large ratio $\mu_{\mathrm{u}}/\mu_{\mathrm{s}} = 10$, say, the magnetic substrate has a destabilizing effect, which is in agreement with earlier findings. For stiffness values in between, the influence changes from stabilizing to destabilizing.

In Fig. 11 we again report the critical stretch (upper panel) and the critical wavenumber (lower panel) for the layered magnetoelastic half-space with $\alpha = \beta = 0.5$, $\mu_{\mathrm{u}}/\mu_{\mathrm{s}} = 0.5, 1, 2, 5, 10$ and the selected values $\bar{B}_{\mathrm{L2}} = 0, 0.5, 1, 1.5, 2, 2.5, 3$. As expected, with no magnetic field the critical values coincide with those shown in Fig. 7. For a stiff upper layer with $\mu_{\mathrm{u}}/\mu_{\mathrm{s}} = 10$, say, the effect of an increase in the magnetic field is minor. For a low stiffness ratio the magnetic field first stabilizes the configuration and, for example, for $\mu_u/\mu_s = 0.5$ an increasing generates instability even in tension.

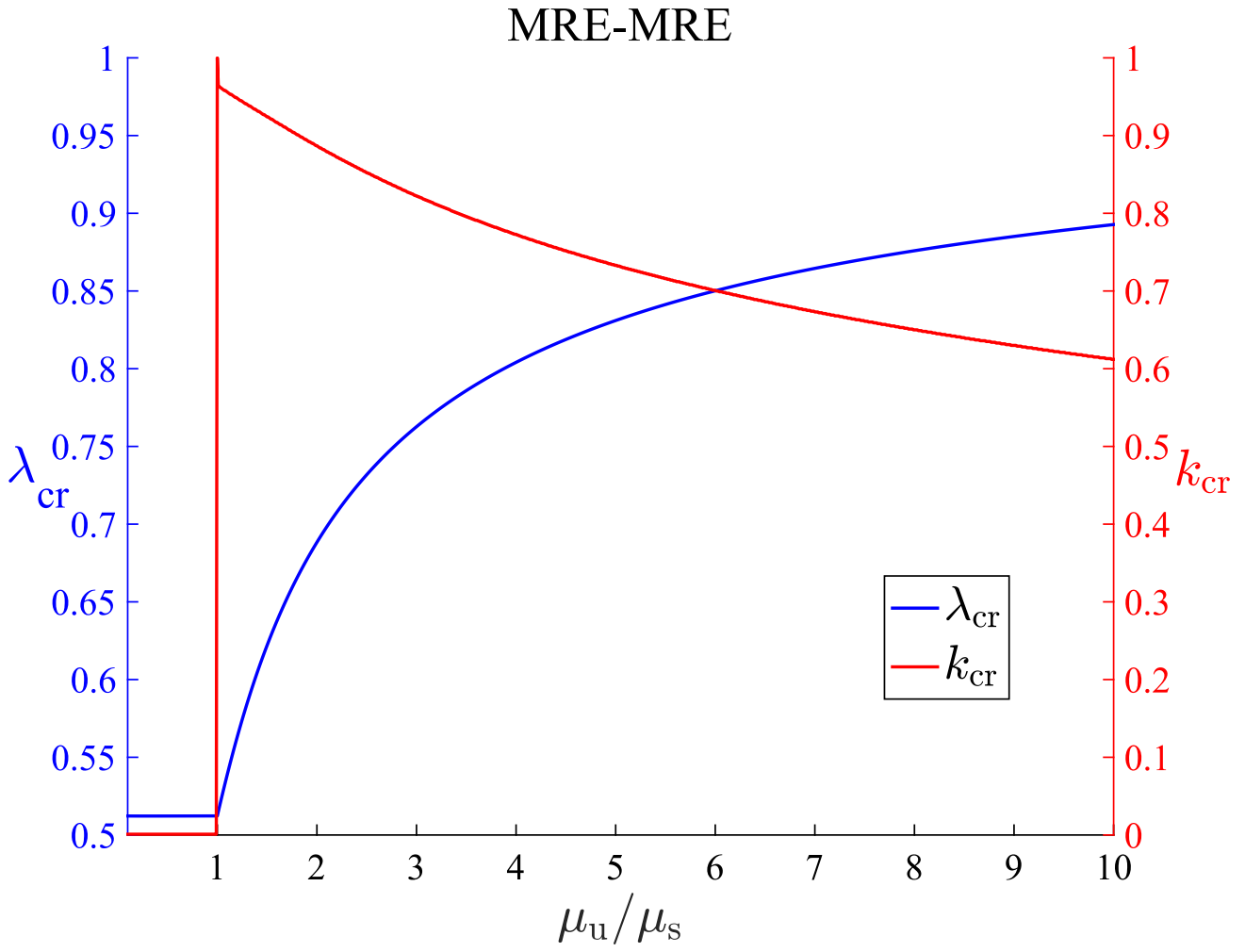

**Fig. 10.** The plot shows the critical stretch $\lambda_{cr}$ on the left and the critical wave number $k_{cr}$ on the right as a function of the stiffness ratio $\mu_u/\mu_s$ for a two-layers magnetoelastic half-space. The results correspond to an applied field $\bar{B}_{L2} = 0.5$ and to the material constants $\alpha = \beta = 0.5$.

The critical wave numbers vanish when $\mu_u/\mu_s < 1$ and the applied fields is small, but increase with the magnetic field, see insert in the lower panel in Fig. 11. For larger ratios, the opposite is true, i.e. the wavenumber gets smaller with increasing field.

In Fig. 12 we evaluate the impact of the magnetoelastic coupling parameters on the critical stretch as a function of the wave number for $\bar{B}_{L2} = 0.5$ and $\mu_u/\mu_s = 5$ (upper panel) and $\bar{B}_{L2} = 2$ and $\mu_u/\mu_s = 5$ (lower panel) for $\alpha = 0.5$ and $\beta = 0, 0.5, 1, 2, 5$. Compared to the results shown in Fig. 9 we again note that the selection of the magnetic properties provides a means to design actuators and controllable devices.

### 8.5. Connections with prior results

Results in Fig. 7 indicate that for a sufficiently stiff magnetoelastic upper layer interacting with softer substrates, an external field $\bar{B}_{L2} \approx 1.5$ may generate surface instability when $\lambda \approx 1$. In Danas et al. (2012) and Psarra et al. (2017), it is shown that wrinkle formation occurs in a ferromagnetic film bonded on a compliant passive foundation for magnetic fields between 0.15 – 0.7 T. The critical magnetic field for commercial magneto-sensitive elastomers with $\mu = 0.1$ MPa, for example, is $B_{L2} = 0.35 - 0.55$ T, which converted to dimensionless form is $\bar{B}_{L2} = 0.98 - 1.55$, hence covered by the range of values considered.

Figure 5b in Psarra et al. (2017) shows data and numerical results of the critical magnetic field as a function of the pre-compression stretch, generating first two-dimensional followed by one-dimensional wrinkling pattern. Periodic wrinkling changes to period-doubling localization for further compression. We were interested in comparing the results obtained from the constitutive model (97) with the experimental and numerical values reported in Psarra et al. (2017). For the upper layer we take $H = 0.8$ mm, $\mu_u = 10$ kPa, for the substrate we use $\mu_s = 3$ kPa. For the magnetic part of the energy function (97) we take $\alpha = \beta = 0.5$. The results, shown in Fig. 13 as a continuous curve, are compared to the experimental data (squares) and the earlier results that account for magnetic saturation (circles). Fair agreement is obtained.

## 9. Concluding remarks

In this article we investigate the surface stability of a layered elastic half-space undergoing large deformation in the presence of a magnetic field normal to its boundary surface. The stability analysis is reduced to evaluating the effect of small increments in deformation and magnetic induction, superimposed on a finitely deformed configuration.

We use the first variation of an energy functional to obtain the equilibrium and magnetic field equations, along with the corresponding boundary conditions. The second variation gives the incremental forms of the total stress, magnetic field, governing equations, and boundary conditions. These are given in Lagrangian form and serve as the basis for the finite element formulation.

For illustration we first consider a layered elastic half-space made of non-magnetizable materials with different stiffness ratios. Then, we analyze a magneto-sensitive half-space subjected to a magnetic induction oriented normal to the free surface and determine the critical stretch for different magnetic properties. We also illustrate the responses of an upper layer resting on a non-magnetizable material and of a two-layer magnetoelastic solid for different mechanical stiffness properties. In each case the numerical results are presented in graphical form.

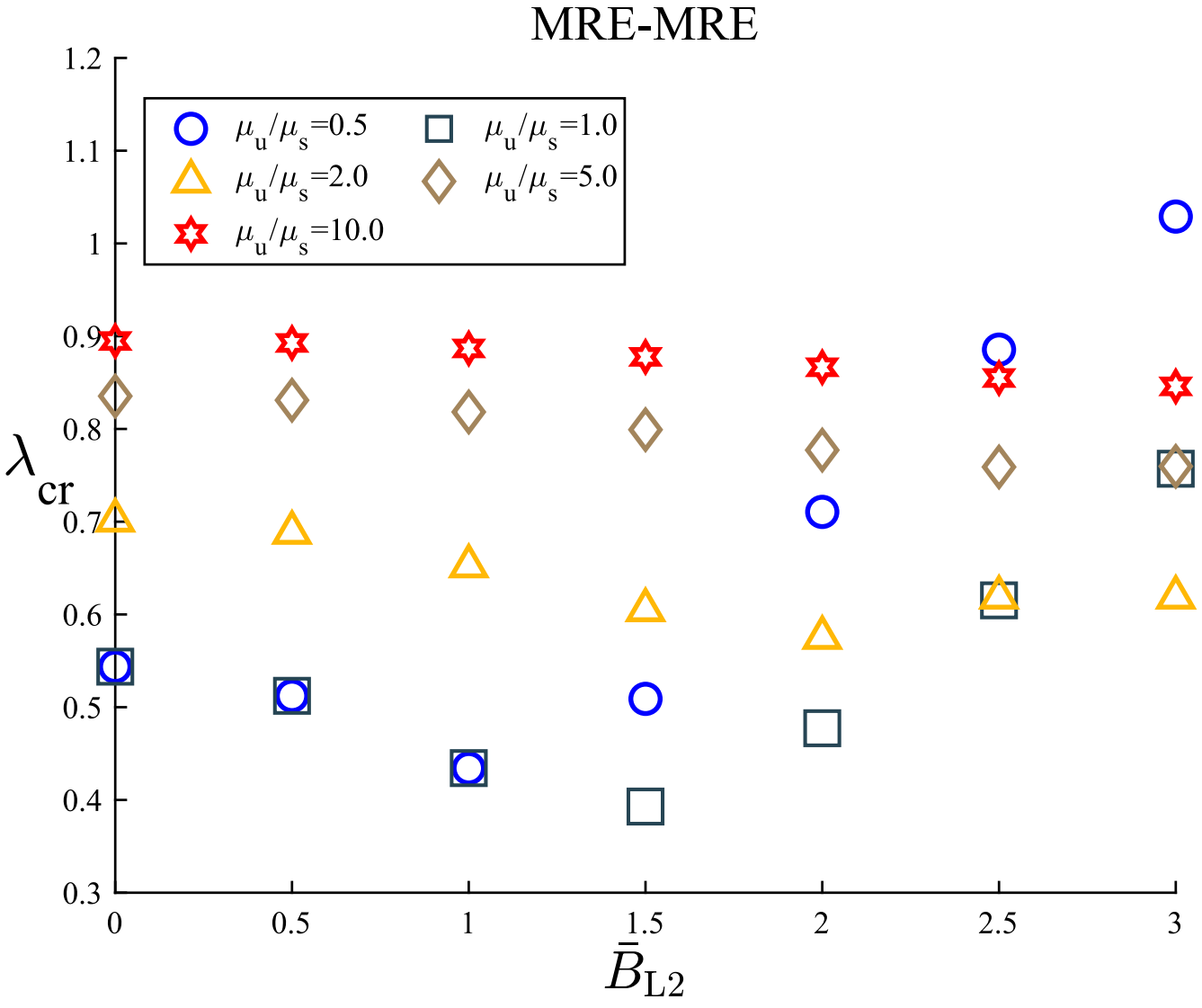

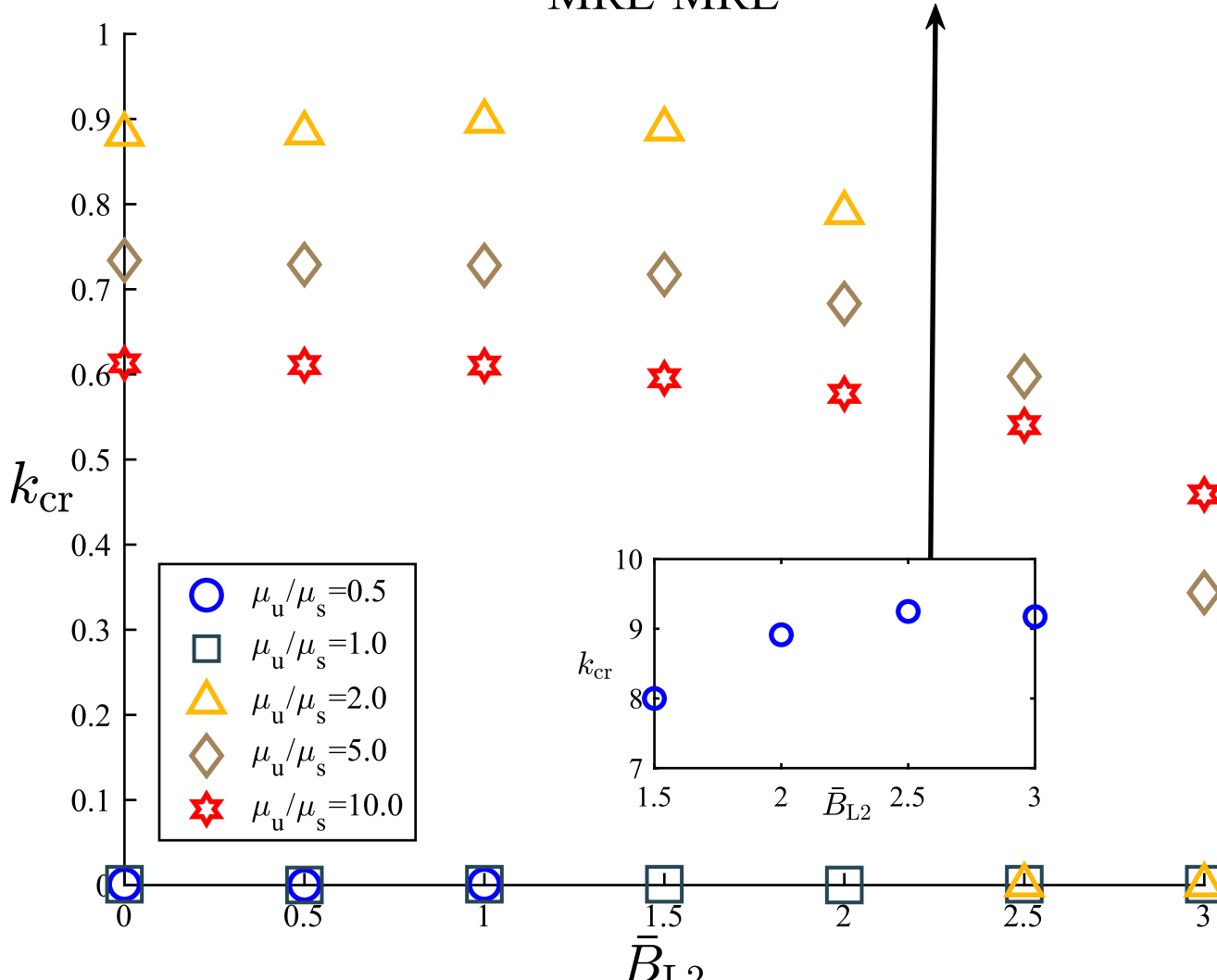

**Fig. 11.** The critical stretch and wave number of a two-layer magnetoelastic half-space for selected values of $\bar{B}_{L2}$ with $\alpha = \beta = 0.5$ and $\mu_u/\mu_s = 0.5, 1, 2, 5, 10$.

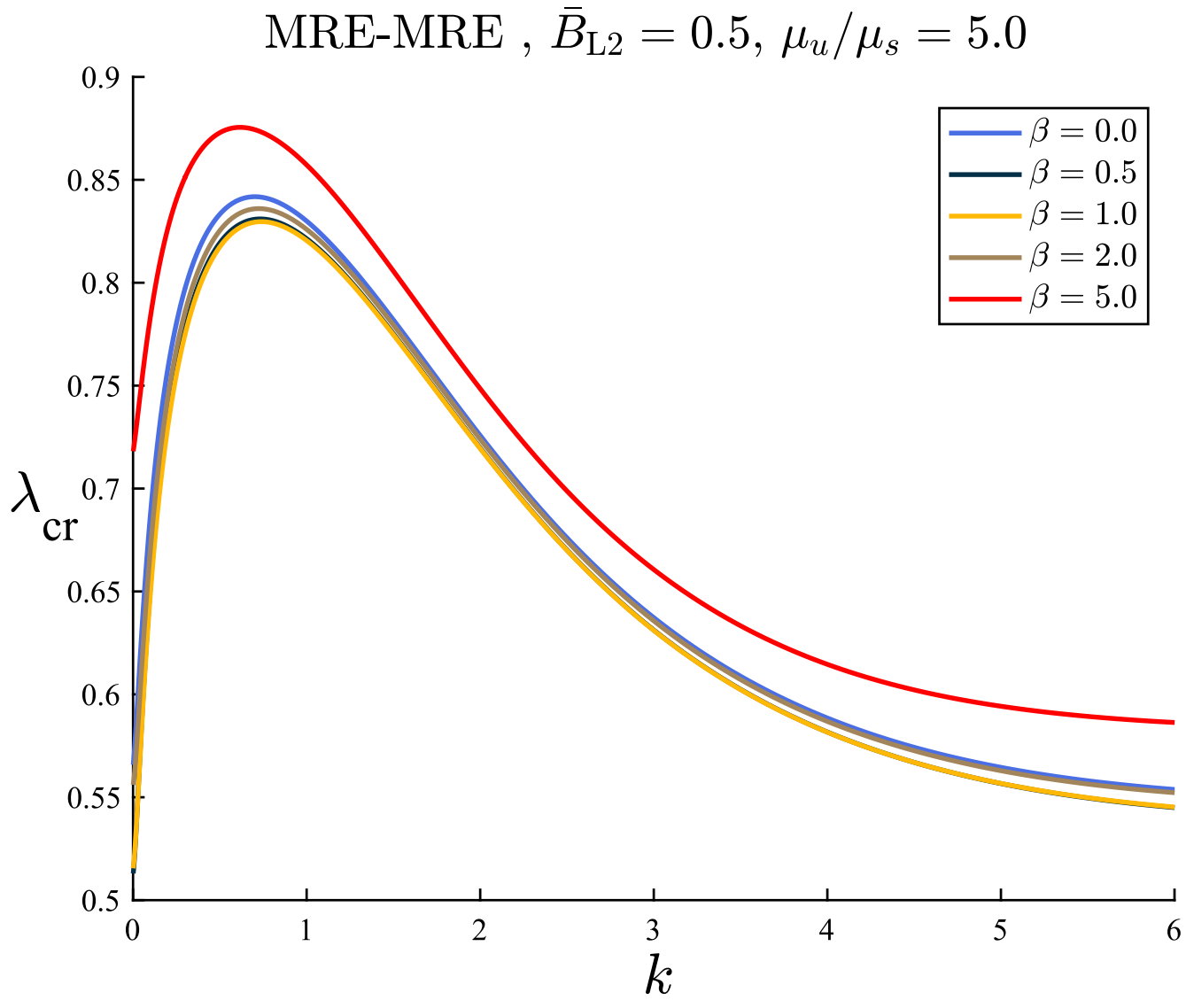

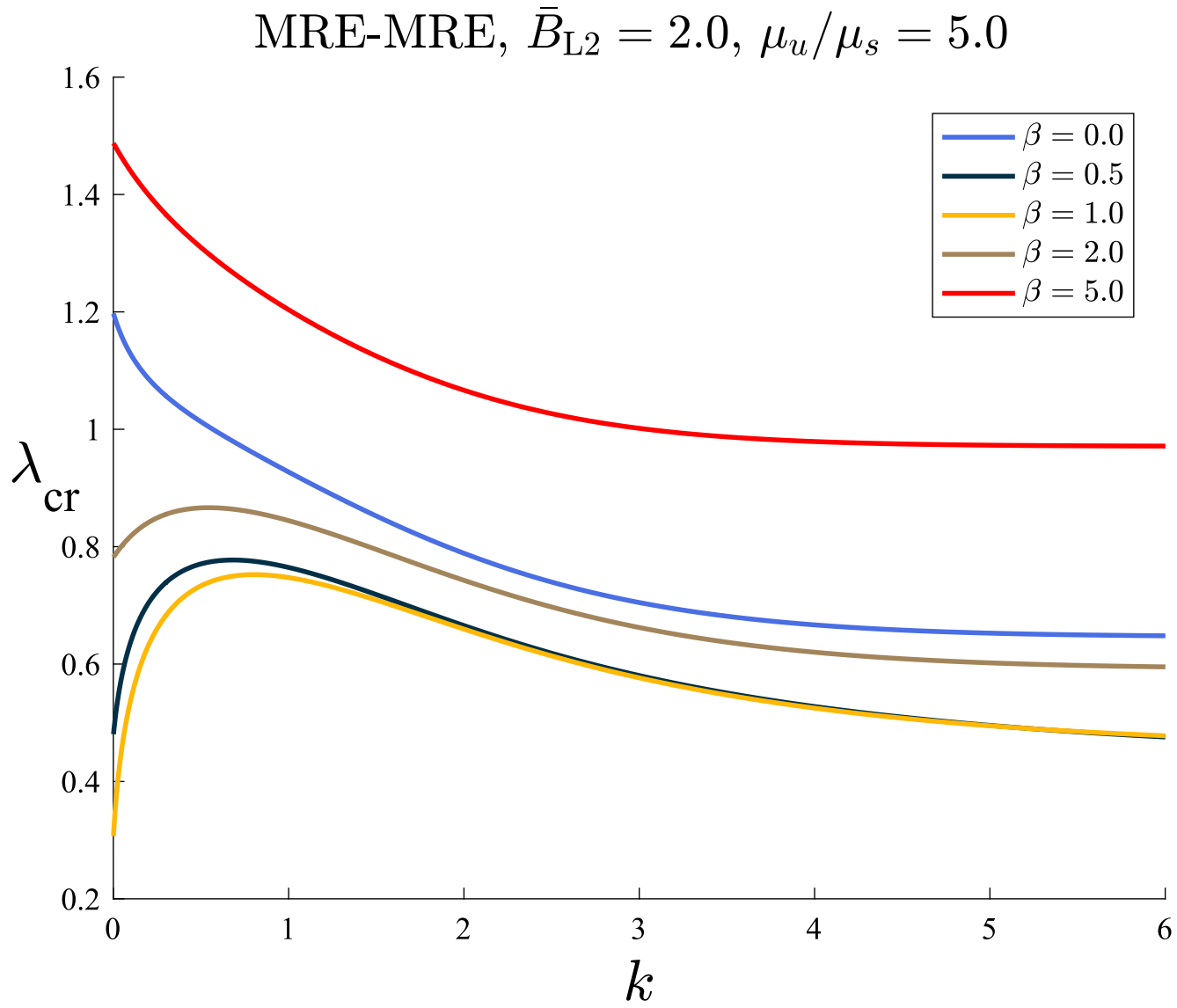

**Fig. 12.** The impact of the magnetic properties on the critical stretch as a function of the wave number for $\bar{B}_{L2} = 0.5$ and $\mu_u/\mu_s = 5$ (upper panel) and $\bar{B}_{L2} = 2$ and $\mu_u/\mu_s = 5$ (lower panel).

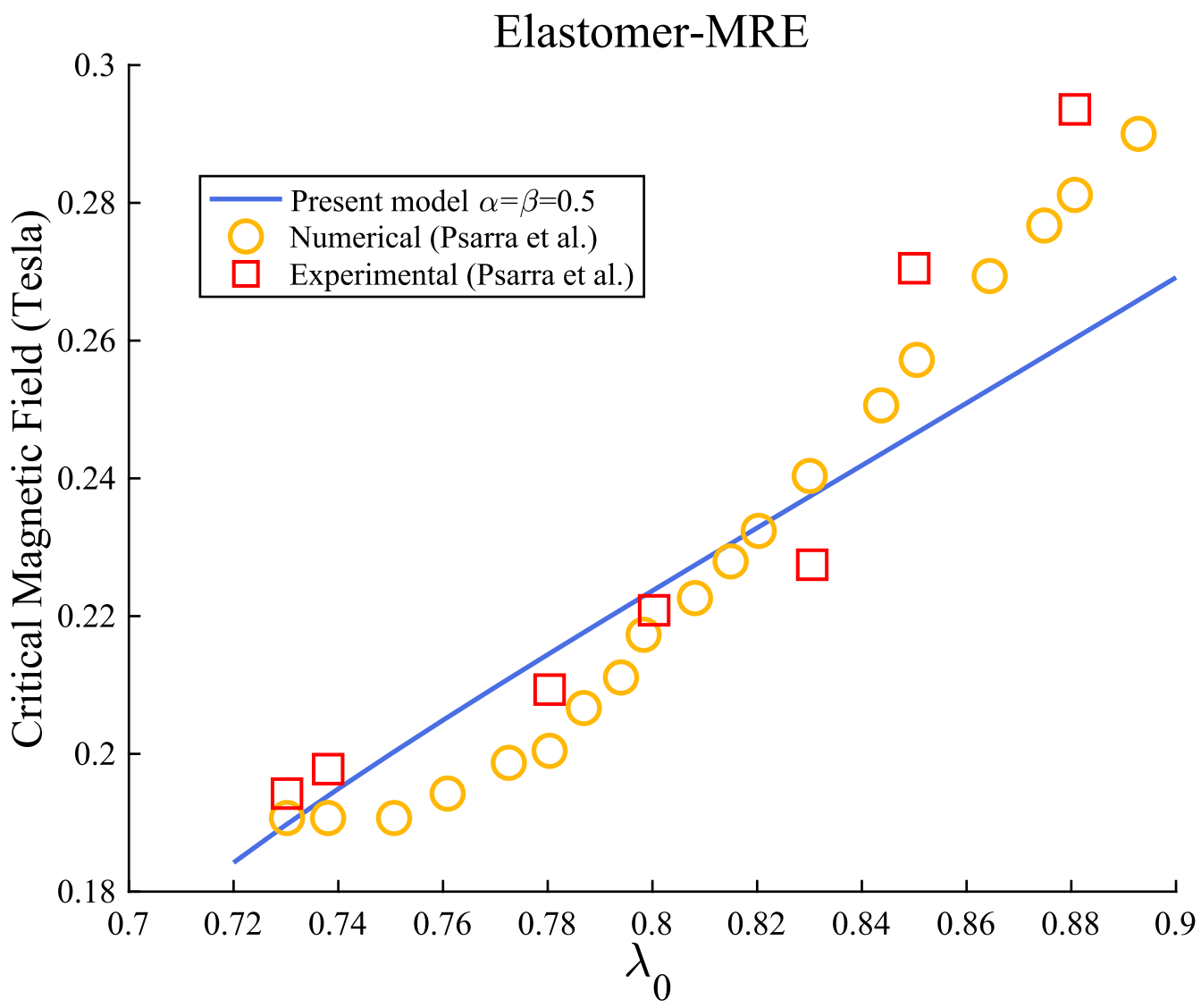

**Fig. 13.** The critical magnetic field as a function of the pre-compression stretch to generate one-dimensional surface wrinkling. We compare results obtained from the constitutive model (97) (continuous curve) with the experimental and numerical values reported in Psarra et al. (2017). Experimental data are shown as squares, the earlier results that account for magnetic saturation are shown as circles.

## CRediT authorship contribution statement

**D. Shahsavari:** Writing – original draft, Visualization, Software, Methodology, Formal analysis. **L. Dorfmann:** Writing – original draft, Supervision, Methodology, Formal analysis. **P. Saxena:** Writing – review & editing, Supervision, Resources, Project administration, Funding acquisition, Formal analysis, Conceptualization.

## Declaration of competing interest

The authors declare that they do not have any conflict of interest.

## Acknowledgments

Davood Shahsavari's work was supported by the Engineering and Physical Sciences Research Council's Doctoral training fund. Davood Shahsavari also acknowledges the Hugh Sutherland and the Jim Gatheral PGR Travel Scholarships to support his visit to Tufts University, USA. Prashant Saxena acknowledges the financial support of the Engineering and Physical Sciences Research Council, UK via project no. EP/V030833/1.

## Data availability

No data was used for the research described in the article.